\documentclass[10pt,letterpaper]{article}
\usepackage[top=0.85in,left=2.75in,footskip=0.75in,marginparwidth=2in]{geometry}

\usepackage[utf8]{inputenc}

\usepackage{cite}

\usepackage{nameref,hyperref}

\usepackage[right]{lineno}

\usepackage{microtype}
\DisableLigatures[f]{encoding = *, family = * }

\raggedright
\usepackage{changepage}

\usepackage[aboveskip=1pt,labelfont=bf,labelsep=period,singlelinecheck=off]{caption}

\makeatletter
\renewcommand{\@biblabel}[1]{\quad#1.}
\makeatother

\usepackage{lastpage,fancyhdr,graphicx}
\usepackage{epstopdf}
\fancyheadoffset[L]{2.25in}
\fancyfootoffset[L]{2.25in}

\usepackage{color}

\definecolor{Gray}{gray}{.25}

\usepackage{graphicx}

\usepackage{sidecap}

\usepackage{wrapfig}
\usepackage[pscoord]{eso-pic}
\usepackage[fulladjust]{marginnote}
\reversemarginpar

\begin{document}
\vspace*{0.35in}

% title goes here:
\begin{flushleft}
{\Large
\textbf\newline{To Infinity and Beyond:
Some ODE and PDE Case Studies}
}
\newline
% authors go here:
\\
P.G. Kevrekidis\textsuperscript{1},
C.I. Siettos\textsuperscript{2},
I.G. Kevrekidis\textsuperscript{3,*},
\\
\bigskip
\bf{1} Department of Mathematics and Statistics, University of Massachusetts Amherst, Amherst, MA 01003-4515, USA
\\
\bf{2} School of Applied Mathematics and Physical Sciences, National Technical University of  Athens, GR 15780, Greece
\\
\bf{3} Department of Chemical and Biological Engineering,
and Program in Applied and Computational Mathematics, Princeton
University, Princeton, NJ 08544 USA; also IAS-TUM, Garching, and
Zuse Institute, FU Berlin, Germany.
\bigskip
* correseponding@author.mail

\end{flushleft}

\section*{Abstract}
When mathematical/computational problems ``reach'' infinity,
extending analysis and/or numerical computation beyond it
becomes a notorious challenge.
We suggest that, upon suitable singular transformations
(that can in principle be computationally detected on the fly)
it becomes possible to ``go beyond infinity" to the other side, with the
solution becoming again well behaved and the computations continuing normally.
%%%
%%%YGK
%%%
%%%  add  citation
%%%
%%%    Wikipedia   https://en.wikipedia.org/wiki/Buzz_Lightyear
%%%
%
In our lumped, Ordinary Differential Equation (ODE) examples this ``infinity crossing''
can happen instantaneously;
%, or may never happen; 
at the spatially distributed, Partial Differential
Equation (PDE) level the crossing of infinity may even persist for finite time, necessitating
the introduction of conceptual (and computational) {\em buffer zones} in which an
appropriate singular transformation is continuously (locally) detected and performed.
These observations (and associated tools) could set the stage for a systematic approach to 
bypassing infinity (and thus going beyond it)  in a broader range of evolution equations;
they also hold the promise of meaningfully and seamlessly performing the relevant computations.
Along the path of our analysis, we present a regularization process via complexification and 
explore its impact on the dynamics; we also discuss a set of compactification transformations and their
intuitive implications..

% now start line numbers
%\linenumbers

% the * after section prevents numbering
\section*{Introduction}
Studying self-similar solutions that collapse in finite time
is a topic of widespread interest in both the mathematical and the physical literature.
The contexts range from scaling~\cite{barenb,Goldenfeld}, to focusing in prototypical dispersive equations such as the
Korteweg-de Vries (KdV) equation~\cite{collapsekdv} and most notably
the nonlinear Schr{\"o}dinger (NLS) equation~\cite{lemesurier,sulem,fibich} on the one hand,
and from droplets in thin films~\cite{thinfilms} and
flow in porous media~\cite{aronson,pelin} to the roughening of crystal surfaces \cite{Kohn} on the other.
One may try to to avoid collapse (e.g. by imposing space~\cite{oberthaler,boris,olshanii}
or time modulations~\cite{saito,porter}) or by identifying the higher order effects
that preclude collapse in physical experiments~\cite{tzortzakis}.
One may alternatively explore what happens to the mathematical, computational
(or even physical~\cite{gaeta}) setting at, or past, the moment of
collapse; see, e.g., the relevant chapter of~\cite{fibich}.
With this latter intent, we start here from an array of simple, controllable
examples and progressively explore more elaborate ones.

Our motivation is simple and, while also physical in part
(as in ~\cite{gaeta}, where the impact
of collapse on optical filaments is sought), it is chiefly mathematical/computational.
As collapse is approached in time, computations naturally break down
and so also do, in part, mathematical approaches; there are notable
exceptions, e.g. efforts to explore beyond collapse, detailed in the
book of~\cite{fibich} for NLS, or in~\cite{aronson,pelin} for the
porous medium problem.
This breakdown has motivated extensive efforts to refine computational
meshes~\cite{renwang,russell} and avoid collapse at the numerical
level (possibly transforming into a co-exploding frame, thus factoring
out the self-similarity~\cite{sulem,siettos} as will be discussed
further below).
Such numerical approaches do not, however, possess the ability to
cross infinity, even in a simpler array of examples
in which we know by construction, or via analytical arguments,
that ``life past infinity persists'' (i.e., that the solution does
not cease to exist and can be continued past a singular point).
This is precisely our aim here: we will propose how to numerically go beyond
infinity as if it was a regular, rather than a singular point.
%(much in the spirit of
%Frobenius methods for series solutions in our undergraduate textbooks).
%
We construct and apply, on demand, a singular transformation that
``absorbs" the singular nature of the dynamics, allowing the solution
to re-emerge on the other side of infinity, where the dynamics
becomes regular again. A complementary perspective
that we explore in this regard is one of compactification transformations which place infinity on equal footing with the rest of the points in, e.g., 
the ODE orbit.

In ordinary differential equation (ODE) examples, an instantaneous encounter with
%/bypassing of 
infinity (crossing or otherwise, as will be discussed below)
will be considered in what follows.
In partial differential equation (PDE) examples, however,
the introduction of physical space leads to multiple possibilities;
one is that collapse might only occur at a single physical point/moment in time, with no
subsequent continuous ``crossing'' of infinity.
This is the so-called transient blowup in the insightful 
summary of~\cite{galakt} aiming at the classification (see the discussion of item (5) therein)
of post-focusing regimes; we will return to it in our discussion. 
A number of important examples have this structure, 
that will not be examined in detail here (although we expect
that the techniques proposed below are quite relevant in these cases too).
Instead, we will focus on the computationally intriguing
case where,  upon "touching" infinity at an initial point in space/time,
the solution will start gradually crossing; in one dimension this will generically result in
two simultaneous crossings that emerge from the original encounter with infinity, and
subsequently propagate apart in space/time.
This poses computational challenges, as collapse persists
in time (there needs to be a singular transformation in some portion(s) of the domain for entire time intervals),
and it is also mobile; we need to adaptively follow
the region(s) where the singular transformation is detected
and accordingly performed as needed. It does not escape us that an additional
possibility can be envisaged: finite {\em spatial} intervals of the solution 
(possibly multiple ones 
simultaneously) may become infinite, leading the regular part
of the solution to be supported in compact regions, resembling
so-called compacton structures originally introduced
in~\cite{rosenau}. This is referred
to as ``incomplete blowup'' in~\cite{galakt}.
Remarkably, complexification of the evolving variable(s)
and how this may lead to a {\em regularization} of the (real) collapsing dynamics
naturally emerges as a potential ``variant'' of collapse.
This is, arguably, a topic of interest in its own merit; yet it connects
naturally with the overall picture of approaching (and potentially crossing) infinity,
and, as such, we will discuss it here.

Our presentation is structured as follows. In
Section II, we will briefly outline our ODE examples and their
concomitant singular transformations, as well as issues of
numerical computation and the notions
of compactification, and of regularization via complexification.
In section III we will discuss one of the
scenarios mentioned above in the case of a 1D PDE.
For this purpose, we "engineer" a transformation of a simple, 1D linear reaction-diffusion problem
that exhibits the ``single initial  $\rightarrow$ multiple mobile" collapse point(s)
scenario; we discuss and illustrate how to address that numerically.
We show in this case too how complexification may lead to regularization.
We then summarize our findings and present some conclusions and topics for future study. Some relevant auxiliary notions are discussed in the 
Appendix and the Supporting Information (SI).

\section*{The ODE Context}

The standard textbook 
ODE for collapse in finite time (and its solution  by direct integration) reads:
\begin{eqnarray}
\dot{x}=x^2 \Rightarrow x(t)=\frac{1}{t^{\star}-t}.
\label{eqn2}
\end{eqnarray}
The collapse time $t^{\star}=1/x(0)$, is fully determined
by the initial condition, and the textbook
presentation usually stops here.
A numerical solver would overflow close to (but before reaching) $t^{\star}$;
yet we can bypass this infinity by appropriately transforming
the dependent variable $x$ near the singularity.
Indeed, the ``good'' quantity $ y \equiv \frac{1}{x} \equiv x^{-1}$,
satisfies the ``good'' differential equation
$\frac{dy}{dt}=-1$;
this equation will help  ``cross'' the infinity (for $x$) by crossing zero
and smoothly emerging on the other side (for $y$).
Once infinity is crossed, we can revert to
integrating the initial (``bad'', but now tame again) equation for $x$.

To manifest the feature that infinity crossing should be thought of as being
on equal footing with any other point on the rest of this orbit, we introduce a notion of
compactification~\cite{compactification}.
Reshuffling the (hyperbolic form of the) solution, we have
\begin{eqnarray}
(t^{\star}-t) x = 1 \Rightarrow
\left(\frac{t^{\star}-t + x}{2}\right)^2 -
\left(\frac{t^{\star}-t - x}{2}\right)^2 =1.
\label{eqn3}
\end{eqnarray}
Compactification through the variables $X$ and $Y$
\begin{eqnarray}
 X &=& \cos(\theta)=(t^{\star}-t-x)/(t^{\star}-t+x),
\label{eqn4}
%\quad
\\
Y &=& \sin(\theta)=2/(t^{\star}-t+x).
\label{eqn5}
\end{eqnarray}
converts this hyperbola to a circle; one can verify
that indeed $-1 \leq X \leq 1$ and $-1 \leq Y \leq 1$ and
also that $X^2+Y^2=1$.
This puts both relevant infinities
\begin{eqnarray}
(t \rightarrow t^{\star}, x \rightarrow
\pm \infty) &\Rightarrow& (X\rightarrow -1, Y \rightarrow 0^{\pm})
\label{eqn6}
\\
(t \rightarrow \pm \infty, x \rightarrow
\pm 0) &\Rightarrow& (X\rightarrow 1, Y \rightarrow 0^{\mp}).
\label{eqn7}
\end{eqnarray}
on equal footing with all other points of the orbit along the
circle.
The trajectory between the point $(1,0)$ (the infinity in $t$, the steady 
state in $x$) and the point $(-1,0)$ (the infinity in $x$) can
be thought of as reminiscent of a "heteroclinic connection".
%
%One might consider this heteroclinic connection between the two infinities as reminiscent of the persistent heteroclinic 
Such connections often arise 
in dynamical systems with symmetries (see e.g. \cite{BackintheSaddle,Scholarpedia}).
The compactification also suggests that, provided we utilize
"the right variables", i.e., the right quantities to observe the solution,
(e.g., in the form of this circle) we should obtain a consistent, smooth picture
(with consistent, smooth numerics).
Indeed, $y(t)=1/x(t) (= t^{\star}-t)$ is a transformation
in itself singular, yet one which converts the ``bad'' exploding variable $x(t)$
into a ``good'' variable $y(t)$, satisfying $\frac{dy}{dt}=-1$
that merely smoothly crosses $0$.

\begin{figure*}[tb]
\begin{tabular}{ccc}
%\begin{center}
\includegraphics[scale=0.18]{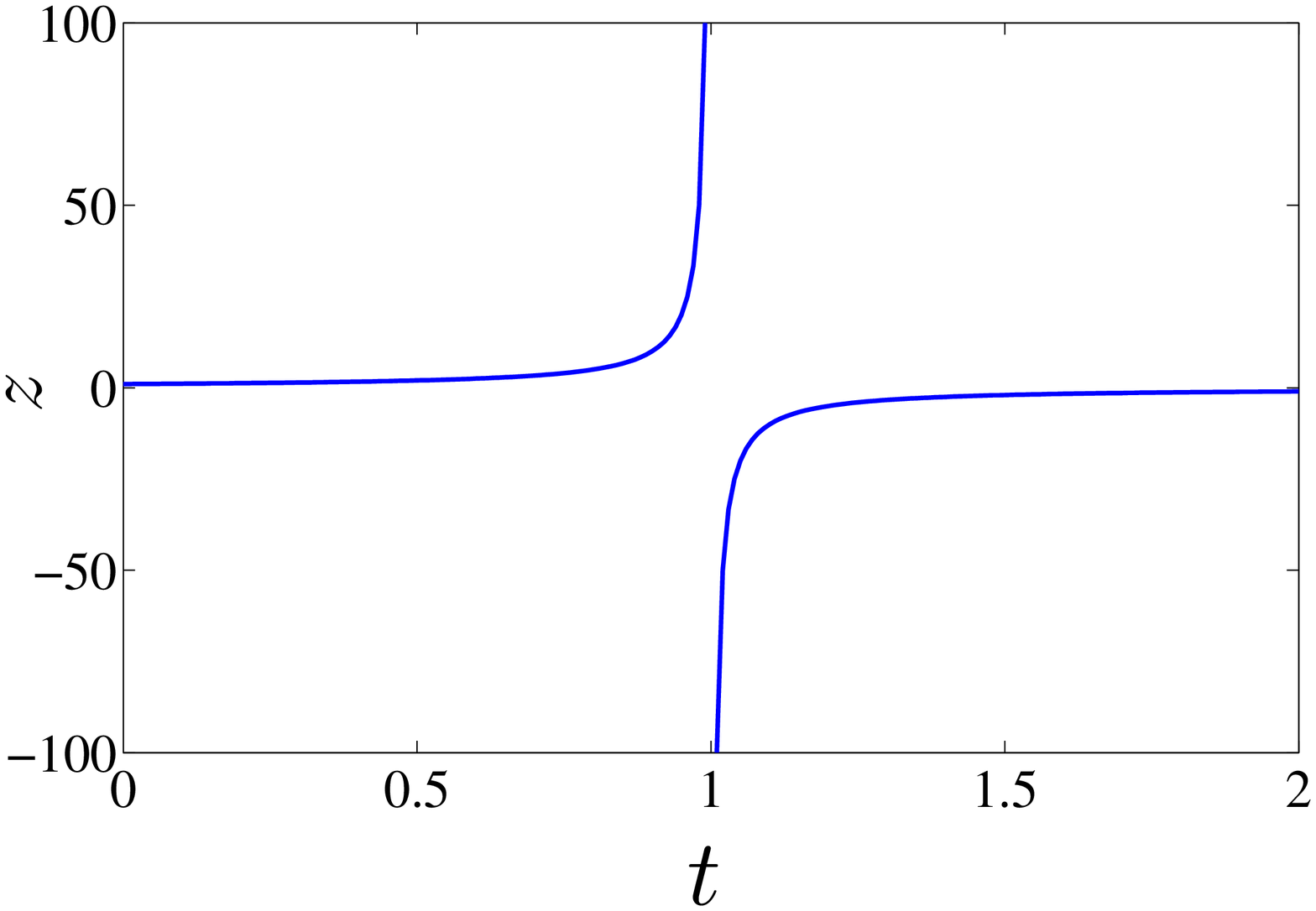} & \includegraphics[scale=0.18]{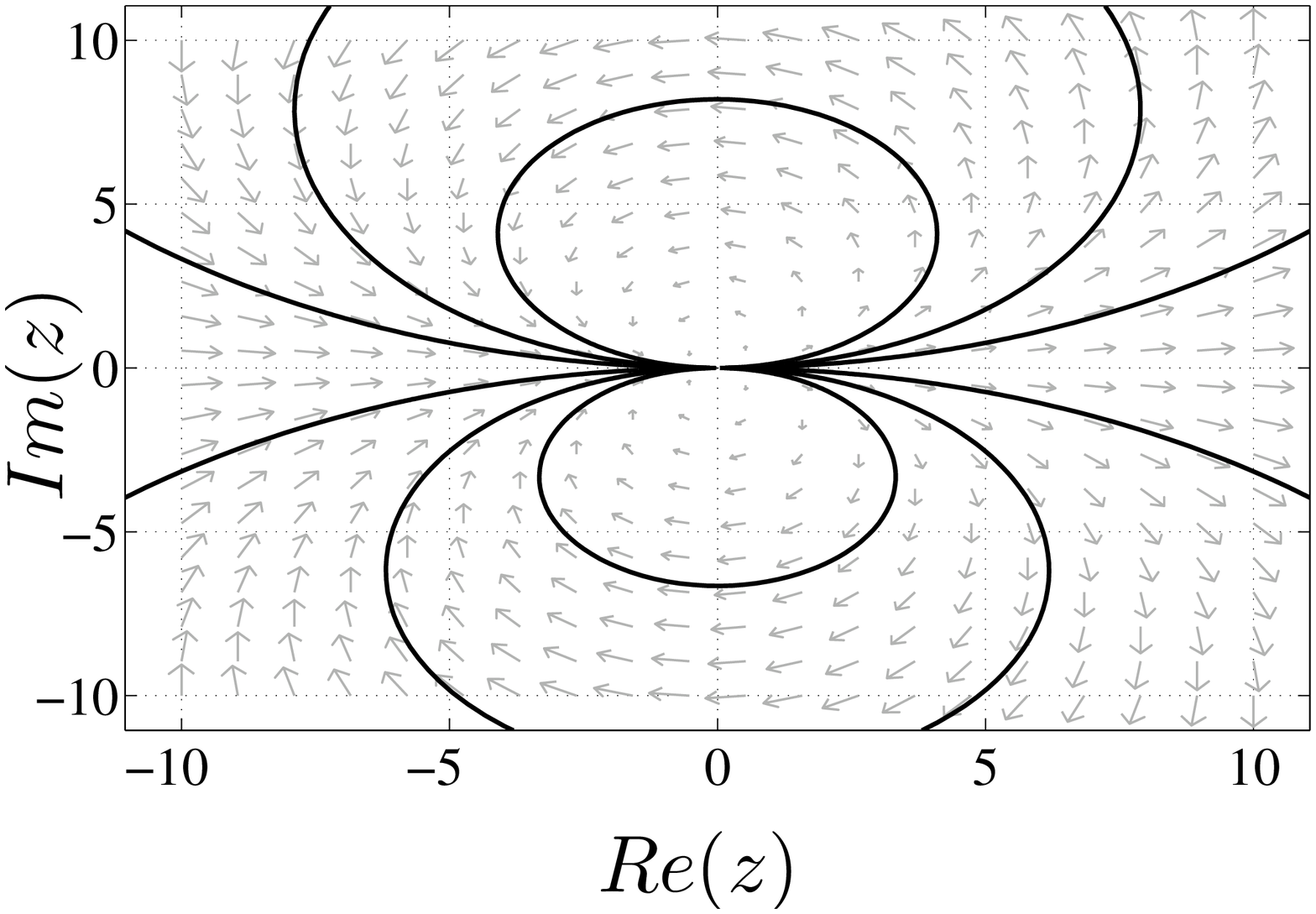} & \includegraphics[scale=0.15]{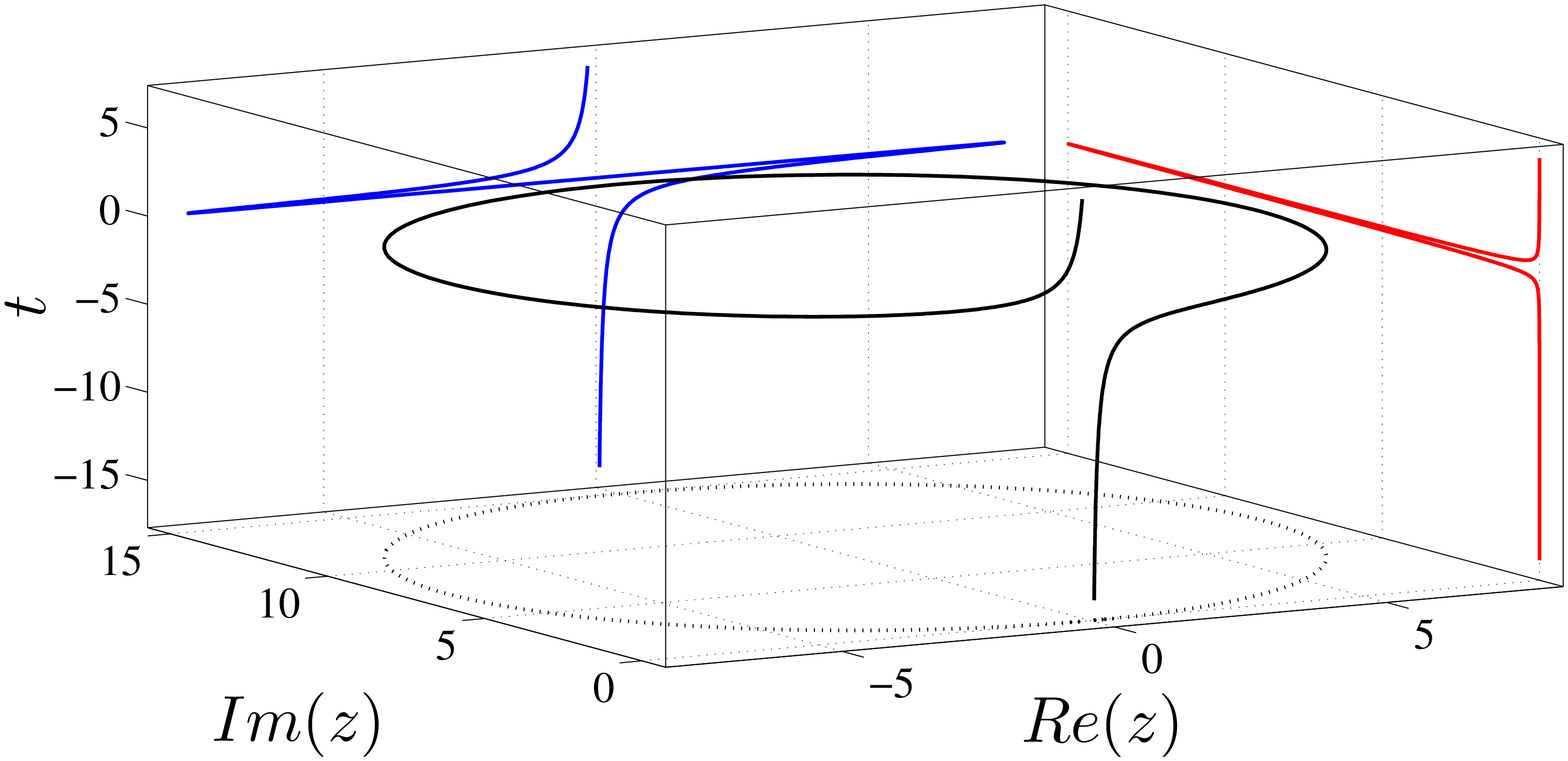}\\
(a) & (b) & (c) \\[8pt]
\includegraphics[scale=0.18]{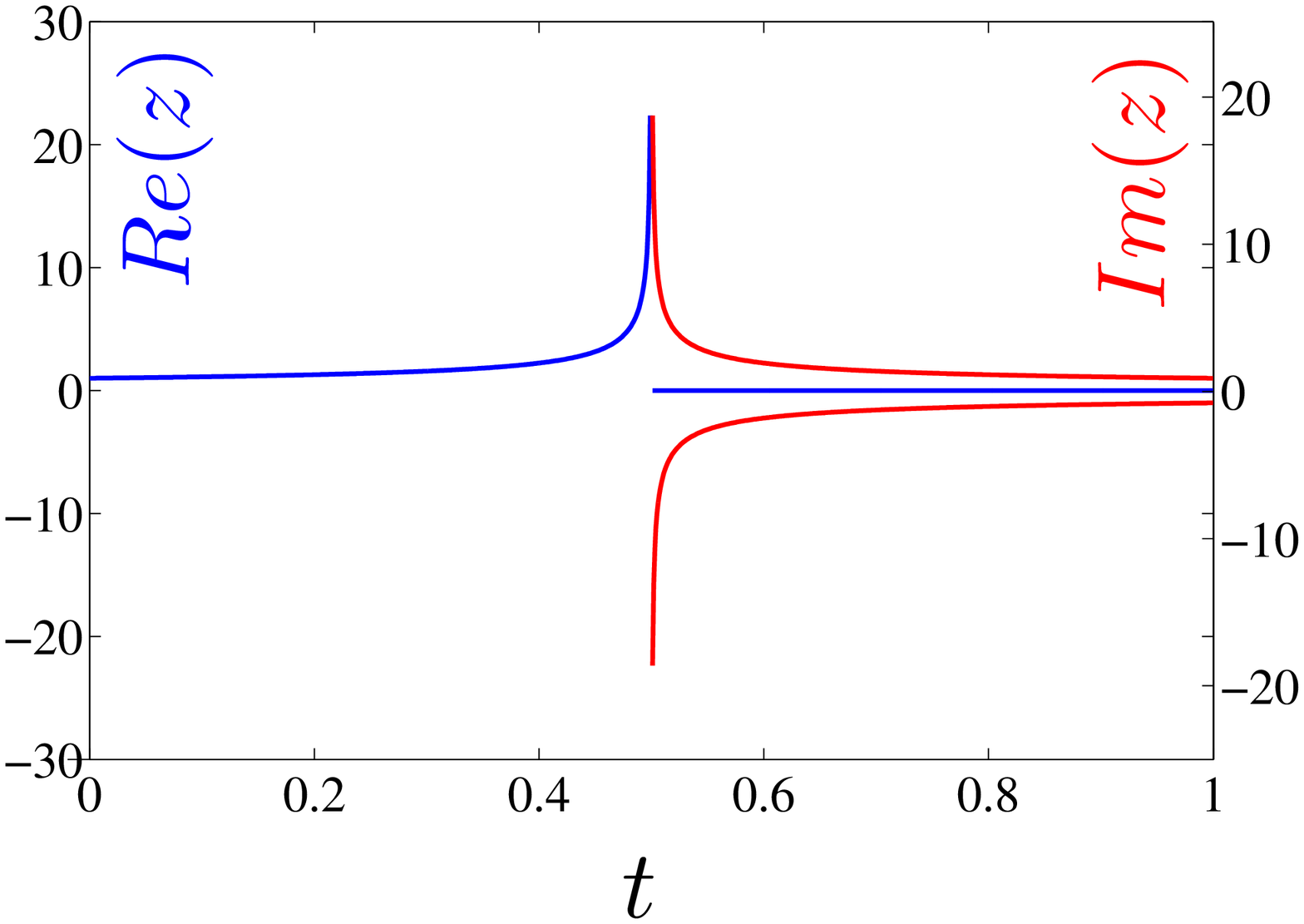}& \includegraphics[scale=0.18]{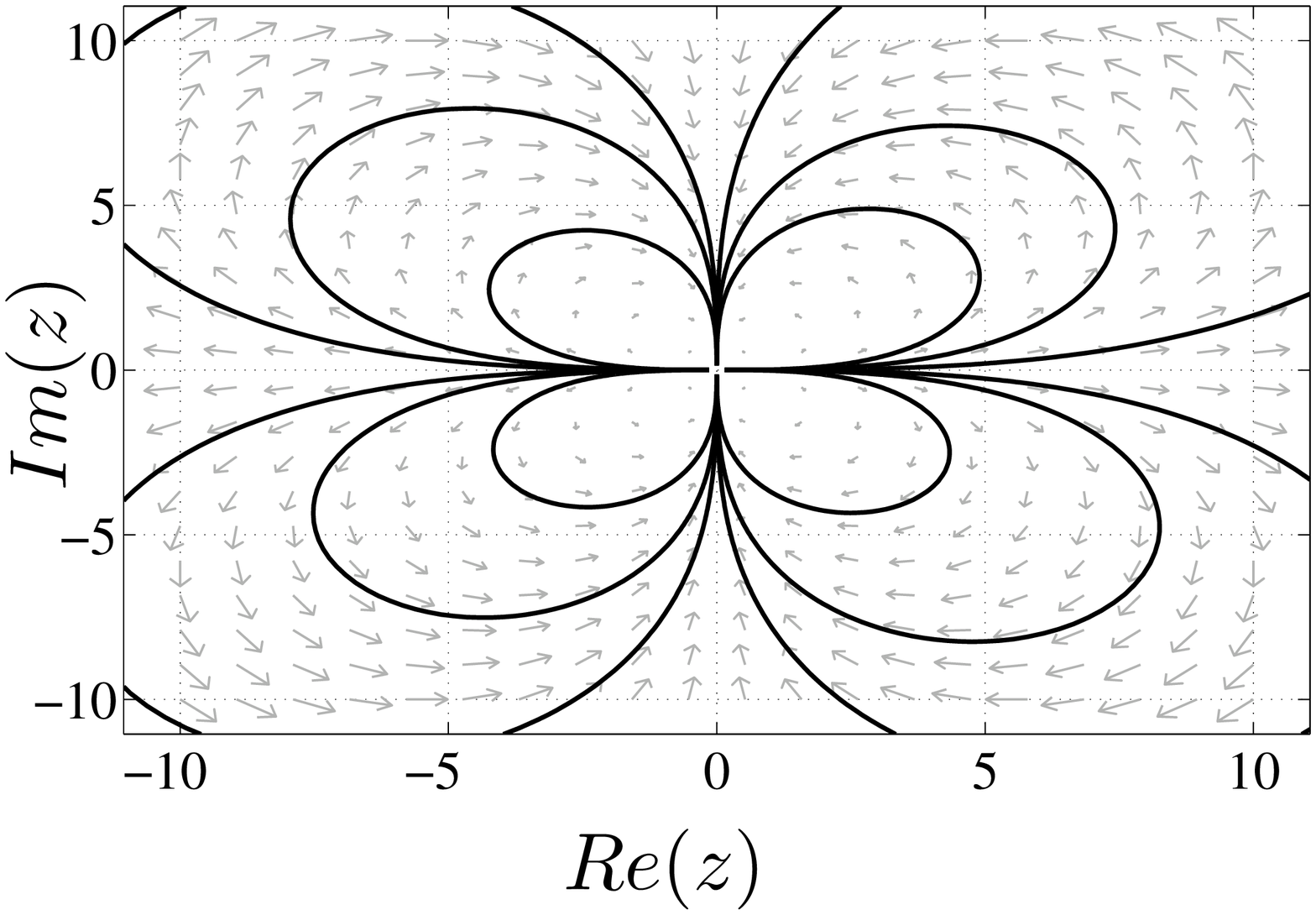} & \includegraphics[scale=0.15]{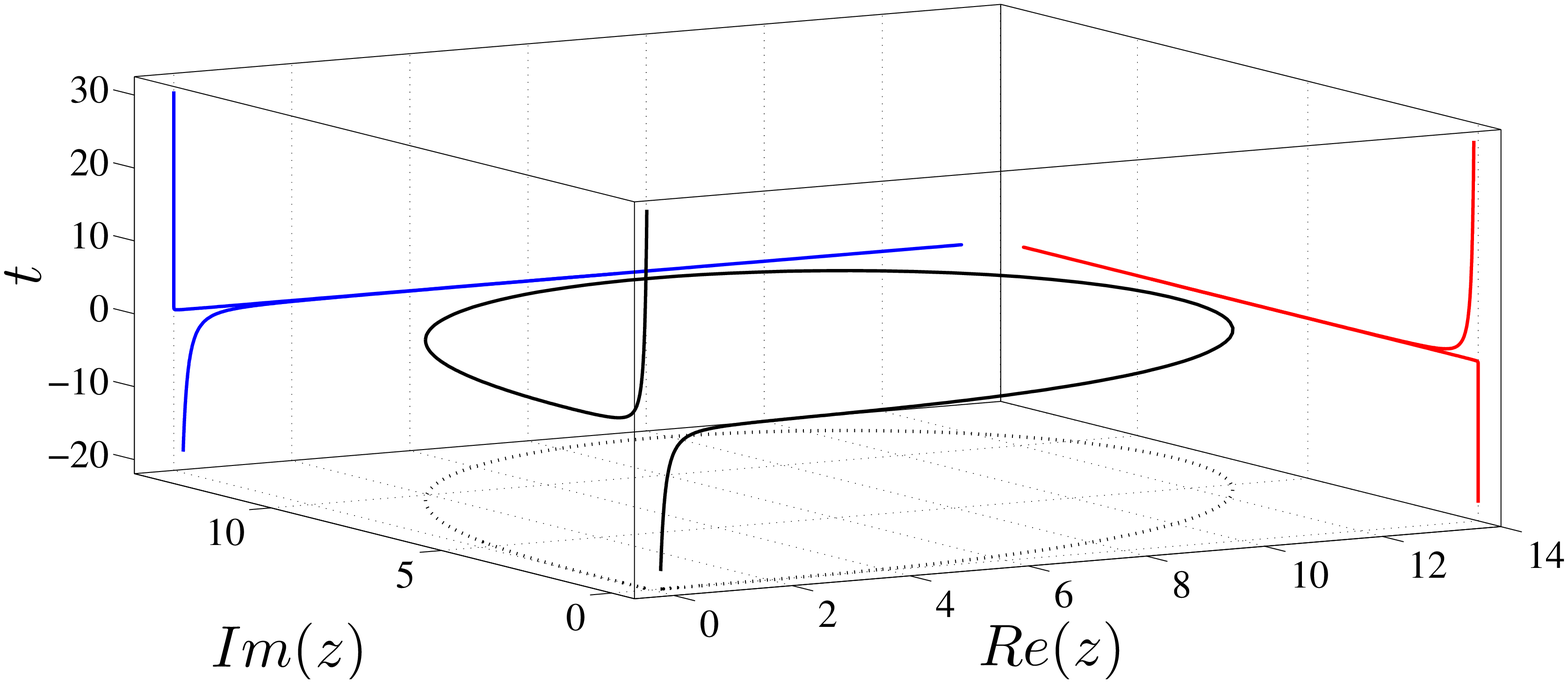}\\
(d) & (e) & (f) \\[8pt]
\end{tabular}
\caption{(Color online) (a) As discussed in the text,
we solve $\dot{x}=x^2$ until the solution
reaches $x(t)=100$, followed by solving Eq.~(\ref{eqn12}) beyond
crossing zero to $y(t)=0.01$, and then returning to Eq.~(\ref{eqn2}).
Here $x(0)=1$, leading to collapse at $t^{\star}=1$. (b), (c) The complex
dynamics of $\dot{z}=z^2$ as represented by ODEs of Eq.~(\ref{eqn21}). Phase plane analysis (b) and sample trajectory (c). Panel (b) shows that the orbits close (and are, in fact, circles as shown in the text). Panel (c) illustrates the circular nature of the projection in the bottom $x-y$ plane, as well as the $x-t$ and $y-t$ plane projections while following the $x-y-t$ composite trajectory.
(d) In this case, we solve $\dot{x}=x^3$ until the solution
reaches $x(t) \approx 25$, followed by solving Eq.~(\ref{eqn12}) beyond
crossing zero to $y(t)=0.04$, and then returning to Eq.~(\ref{eqn2}).
Here $x(0)=1$, leading to collapse at $t^{\star}=0.5$. 
%(b), (c) The complex
%dynamics of $\dot{x}=x^3$ as represented by ODEs of Eq.~(\ref{eqn21})
(e), (f) The complex dynamics of the
two-degree-of-freedom system with $\dot{z}=z^3$.
An example from the first quadrant of panel (e) is illustrated in more detail in the panel (f), exhibiting how collapse is avoided in this case.
}
\label{fig2}
%\end{center}
\end{figure*}

The following numerical protocol then naturally
circumvents problems associated with infinity
in a broad class of ODEs that collapse self-similarly,
as power laws of time (or, importantly, as we will see below including
in the SI, also
{\em asymptotically} self-similarly):

\begin{itemize}
\item Solve the ``bad'' ODE of Eq.~(\ref{eqn2}) for a while,
continuously monitoring, during the integration, its growth towards collapse.
\item If/when the approach to collapse is detected, estimate its (asymptotically)
self-similar rate (the exponent of the associated power law, here $-1$) and use it to
{\em switch} to a ``good'' equation
for $y$, relying on the singular transformation $y=1/x$ with this exponent
(and on continuity, to obtain appropriate initial data for
this good equation).
\item 
Run this ``good'' equation for $y$ until $0$ for it (or $\infty$ for the former, ``bad'' equation) is
safely crossed, computationally observing for $x$ an (asymptotically) self-similar ``return'' from infinity.
\item Finally, transform back to the ``bad'' equation (no longer that bad, as infinity
has been crossed) and march it further forward in time.
\end{itemize}
This protocol has been carried out in Fig.~\ref{fig2} (see caption), illustrating that
the dynamics can cross infinity and computation can be continued
for all time, provided that the self-similar approach to infinity is adaptively detected and
the associated, and appropriately numerically estimated, singular transformation is used to cross it.
Note that the compactification also  
allows the progression past infinity {\em in time} too, when now $y$ crosses
zero as time approaches positive infinity and then ``returns'' from negative
infinity.
%
%The approach allows a type of "flow box" transformation to be brought to bear on the entire vector field, {\em including
%the steady state at $x=0$ and the ``point at infinity''}: a constant %%(transformed) time evolution on the constant curvature circle. 

%\begin{figure}[bt]
%\begin{center}
%\vspace{-0.1cm}
%\includegraphics[height=.22\textheight, angle =0]{stathis1.eps}
%\end{center}
%\par
%\vspace{-0.7cm}
%\caption{Solve Eq.~(\ref{eqn2}) until the solution
%reaches $x(t)=100$, then solving Eq.~(\ref{eqn12}) beyond
%crossing zero to $y(t)=0.01$, and then returning to Eq.~(\ref{eqn2}).
%Here $x(0)=1$, leading to collapse at $t^{\star}=1$.}
%\label{fig1}
%\end{figure}
%

%
We can now try to extend/generalize these ideas to other collapse rates
(i.e. arbitrary powers/exponents of self-similarity).
The collapse of $\dot{x}=x^3$  whose exact solution is $x(t)=\frac{1}{\sqrt{t^{\star}-t}}$
is worth examining separately.
The relevant singular transformation
(here $y(t)=1/x^2$) will take us to infinity in finite time,
but, at first sight, will not cross - $y(t)$ becomes imaginary
beyond $t^{\star}$.
An appropriate compactificaton resolves the issue
\begin{eqnarray}
 X &=& \cos(\theta)=(t^{\star}-t-x^2)/(t^{\star}-t+x^2)
\label{eqn11}
\\
Y &=& \sin(\theta)=2/(t^{\star}-t+x^2),
\label{eqn12}
\end{eqnarray}
leading to perfectly regular dynamics on a circle,
so that the singularity is again ``bypassed"  in the spirit
of Fig.~\ref{fig2}.
%%%%%%%%%%%%%%%%%%%%%%%%%%%%%%%%%%%%%%%%%%%%
%%%
%%% Itako, here something funny happens
%%% I worry about the TIME that it take for the imaginary
%%% part to go to zero and "give over" to the real part
%
% PGK: I think this looks fine.
Yet an implicit multi-valuedness clearly arises as a crucial issue
in selecting useful transformations for such exact solutions;
we will return to this important issue below.

For the time being, we argue that one can generalize the above notions
for ODEs that asymptotically collapse self-similarly, $x(t) \sim 1/(t^{\star}-t)^a$, so as to produce a useful compactification in the form
\begin{eqnarray}
 X &=& \cos(\theta)=((t^{\star}-t)^a-x)/((t^{\star}-t)^a+x)
\label{eqn13}
\\
Y &=& \sin(\theta)=2/((t^{\star}-t)^a+x).
\label{eqn14}
\end{eqnarray}
In this form, the dynamics ``benignly'' travels along the circle.
Relevant examples can straightforwardly be extended
to, e.g., fractional powers although it is known from 
standard ODE analysis that issues of uniqueness may arise
there that we do not delve into in the present work.
%
%Alternative useful compactifications, bypassing issues with %radicals
%are discussed in Appendix/Supplementary Information....

%
More generally then, the self-similarly collapsing ODE
$\frac{dx}{dt}= \pm x^p$ has the solution
$\pm \frac{1}{1-p} x^{1-p}= t-t^{\star}$
and its scaling in time follows $x(t) \sim (t^{\star}-t)^{1/(1-p)}$,
with the collapse time once again determined by the initial data.
Given a ``legacy code'' that integrates the ODE $\dot{x}=F(x)$,
we monitor its growth approaching collapse (i.e., how $F(x)$ scales as $x^p$,
or more generally with $||x||$)~\footnote{For vector cases, the analogous feature will be to
monitor
the norm dependence as $||x||^p$, although we will not explore
such a case here.}.
%%%%%%%%%%%%%%
%%%%%%%%%%%
Upon detection of asymptotically self-similar collapse, at sufficiently large $|x|$
(e.g $10^{2}$ in the ODE of Fig.~\ref{fig2}, or $10^4$ in the PDEs of the
next section) we stop
solving the ``bad'' ODE.
We use instead the
singular transformation $y=x^{-p+1}$ (more generally
$y=\int 1/f(x) dx$) to solve the ``good'' $y(t)$
ODE that crosses $0$ rather than infinity.
Then, a little beyond the collapse time (beyond infinity for $x(t)$,
beyond $0$ for $y(t)$) we simply revert to the original,
``bad'' (yet no longer dangerous~!) ODE, with continuity furnishing the 
relevant matching conditions.
An illustration of asymptotically self-similar blowups, where different
transformations are used to cross two different infinities (the finite-time/infinite value and the infinite-time/finite value ones) is included in the SI. 

Examining such infinity crossings
as regular, rather than singular points begs an ``explanation'' for the mechanism of
exiting the real axis along $+\infty$ and then re-emerging on the
other side at $-\infty$  (for $\frac{dx}{dt}= x^2$)
or -arguably more remarkably-  from $+ i \infty$
back towards the origin in the example involving $x^3$).
In the latter case there is an obvious ambiguity: the solution might just as well be chosen to re-emerge from $-i \infty$: one can formally, past the collapse, accept $x(t)=\sqrt{-1/(t-t^{\star})}=i/\sqrt{t-t^{\star}}$ for
$t>t^{\star}$ or, alternatively, $x(t)=\sqrt{1/(-(t-t^{\star}))}=-i/\sqrt{t-t^{\star}}$.
This is perhaps a prototypical (and tangible) example of the ``phase loss''
feature argued in~\cite{gaeta,fibich}.
%
%The compactification approach developed above provides a rationalization of this
%feature: it  ``wraps around'' the real axis (for $x$ in Eq.~(\ref{eqn2}), and
%for $x^2$ in the cubic example, justifying the
%regular (and infinitesimal in time!) crossing "from one infinity to the other".

As a way of shedding further light into these features, we examine the
{\em complexified version} of Eq.~(\ref{eqn2}).
The complexified version $\dot{z}=z^2$ leads to the
two-dimensional dynamical system:
\begin{eqnarray}
\dot{x} = x^2-y^2, \quad         \dot{y} = 2 x y.
\label{eqn21}
\end{eqnarray}
The real axis is an invariant subspace, retrieving our real results;
yet complexification endows the dynamics with an intriguing
``capability'':
as Fig.~\ref{fig2} illustrates through the $(x,y)$ phase plane,
collapse {\em is avoided} in the presence of a minuscule imaginary part.
Large ``elliptical-looking'' trajectories are traced
on the phase plane, eventually returning to the neighborhood of the sole fixed
point of $(0,0)$ --which in the real case one would characterize as semi-stable--.
%%%
%%%YGK
%%%   Itako, what is the 2D eigenvalue ? is it now semistable in the 
%%%   complex case ? Does this mean that there will be longer ad looonger
%%%   intervals for every return, like the Guckenheimer cycles ?
%%%%%%%
% PGK: Eigenvalue is 0. I presume that all returns take 
% infinite time as 0 is approached asymptotically.
%
The system of Eq.~(\ref{eqn21}) can
be tackled in closed form since the ODE $\dot{z}=z^2$ yields $\frac{1}{z}=-t + \frac{1}{z(0)}$.
For $z=x+i y$ ($z(0)=x_0+i y_0$) we obtain the explicit orbit formula
\begin{eqnarray}
x(t) &=& \frac{x_0 (x_0^2+y_0^2)- t (x_0^2+y_0^2)^2}{(x_0 - t (x_0^2
+y_0^2))^2 + y_0^2}
\label{eqn23}
\\
y(t) &=& y_0 \frac{x_0^2 + y_0^2}{(x_0 - t (x_0^2
+y_0^2))^2 + y_0^2}.
\label{eqn24}
\end{eqnarray}
Eliminating time by dividing
the two ODEs within Eq.~(\ref{eqn21}) directly yields an ODE for $y=y(x)$
(rather than the parametric forms of Eqs.~(\ref{eqn23})-(\ref{eqn24})).
From this ODE, one can obtain that the quantity
\begin{eqnarray}
E=\frac{y^2+x^2}{y}=\frac{y_0^2+x_0^2}{y_0}
\label{eqn25}
\end{eqnarray}
is an {invariant} of the phase plane dynamics,
and thus the latter can  be written as $x^2 + (y-R)^2 = R^2$,
where $R^2=(x_0^2+y_0^2)^2/(4 y_0^2)$.
%
%%%YGK  I added a square on the last R above -- pleased check that OK
%%%
%% PGK: Checked.
%
That is, the trajectory evolves along circles of radius $R$
in the upper (resp. lower) half-plane if $y_0 > 0$ (resp. $y_0<0$.)
Approaching the axis
with $y_0 \rightarrow 0$, the curvature of these circles tends
to $0$ and their radius to $\infty$ (retrieving the real dynamics
as a special case).
Fig.~\ref{fig2} through its planar projections
%%%
%%%YGK
%%%  do a nice "word job" about the description 
%%%   with words.  It is now all in the figure caption
%%%  maybe say "it is all explained in the figure caption
%%%  and fix it there --
%%%%%%%%%%%%%%%%%%%%%%%%%%%%%%
illustrates not only the radial projection of the dynamics
in the $x-y$ plane, but the $x-t$ and $y-t$ dependencies.
Starting with a minuscule imaginary part the real
dynamics tends to infinity; yet when the real part gets sufficiently large
(somewhat in the spirit of our computations above),
the imaginary part ``takes over'', grows rapidly, and ``chaperons''
the real part to the negative side.
Once the real part reaches the ``opposite" (absolutely equal)
negative value, the imaginary part rapidly shrinks
and the formerly bad, yet now benign real equation ``takes over'' again.

We point out here that there is also a canonical way to generalize the
compactification of this complex picture to the Riemann sphere
through the inverse stereographic projection
\begin{eqnarray}
X &=& \frac{2 x}{x^2+y^2+1}
\label{eqn26}
\\
Y &=& \frac{2 y}{x^2 + y^2 + 1}
\label{eqn27}
\\
Z &=& \frac{x^2 + y^2 -1}{x^2+ y^2 + 1}.
\label{eqn28}
\end{eqnarray}
Now the real dynamics become a great
geodesic circle, while all other complex plane curves become regular circles
on the surface of the sphere.
Under this transformation all points with $x^2+y^2 \rightarrow \infty$
are identified with $(0,0,1)$, rationalizing
the vanishing time needed to move from one to the other. 

For the cubic case $\dot{z}=z^3$ the two-dimensional dynamical system becomes
\begin{eqnarray}
\dot{x} &=& x^3 - 3 x y^2
\label{eqn29}
\\
\dot{y} &=& 3 x^2 y - y^3.
\label{eqn30}
\end{eqnarray}
%%%https://preview.overleaf.com/public/dcgjpqcxrxnz/images/9d1cfe61623902c54f4fffe0e99ca185a66194ce.jpeg
The corresponding phase portrait is shown in
Fig.~\ref{fig2}(c), while a typical trajectory
is shown in figure Fig.~\ref{fig2}(f).
Instead of one ``leaf'' in the upper half plane there are now two leaves,
one in each quadrant, see Fig.~\ref{fig2}(d); this suggests a natural generalization to $n-1$ leaves in each half plane
in the case of $\dot{z}=z^n$.~\footnote{It does not escape
us here that a particularly intriguing case in its own
right is when $n$ is rational and perhaps even more so when it is
irrational. However, we will restrict our considerations to
the simpler integer cases herein, deferring the rest to
future work.}
In the cubic case there is collapse
for both positive and for negative initial data, and reentry
along either the positive (resp. the negative) imaginary 
infinity (i.e., from $+i \infty$, resp.$-i \infty$)
could be chosen (in analogy to the arbitrariness in the phase factor).
%
%%%
%%%  YGK
%%%   Itako, do you have a sense if this corresponds in some way to
%%%   the Minkowski, de Sitter or anti-de-Sitter geometries ?
%%%   I know this is crap, but there ought to be an analogy in these
%%%   possible compactifications with the possible vacua of cosmology
%%%
%%% PGK: I don't see any connection that I can properly quantify beyond "words"
%%% 
%%%  YGK:  OK - skip for now...
%%%
%%%%%%%%%%%%%%%%
However, for even infinitesimally small imaginary data, the  
%${\mathcal PT}$
symmetry is broken, and unique trajectories are selected along
each quadrant.
A small real part (accompanied  by a small imaginary part)
as in the bottom right panel
of Fig.~\ref{fig2} grows until eventually (when sufficiently large)
the imaginary part takes over.
The real part then decays rapidly to $0$, while the imaginary decays
slowly, closing the orbit in the first quadrant; this is again
a natural extension of the limiting case of purely real initial data.
This complex formulation also allows the quantification 
(in a vein similar to~\cite{sntime,bender}) of how
long it takes for initial data, say, on the real axis,
to ``emerge" on the imaginary axis. 
In the SI we show that this time tends to
$0$ for the transitions between $+\infty$ and $+ i \infty$ for $\dot{z}=z^3$ (or
from $+\infty$ to $-\infty$ in $\dot{z}=z^2$).
%

%It is tempting to attempt a qualitative analogy between the dynamics observed here (excursion to infinity and return to zero)
%and the persistent heteroclinic cycles observed in systems with symmetries, like the
%Kuramoto-Sivashinsky PDE \cite{BackintheSaddle,Guckenholmes,Scholarpedia}.
%
%The compactified infinities at the top of the Riemann sphere play the role of one of the steady states involved,
%while zero, the ``bottom'' of the Riemann sphere, plays the role of the other. Such ``heteroclinic connections" between finite states and infinities arise also in PDE cases as we will discuss below in the NLS case, as well as other problems with symmetry, e.g. %\cite{Knob1,Knob2}.

\section*{A PDE Case}
\begin{figure*}[tb]
\begin{tabular}{cc}
%\begin{center}6
\includegraphics[scale=0.4]{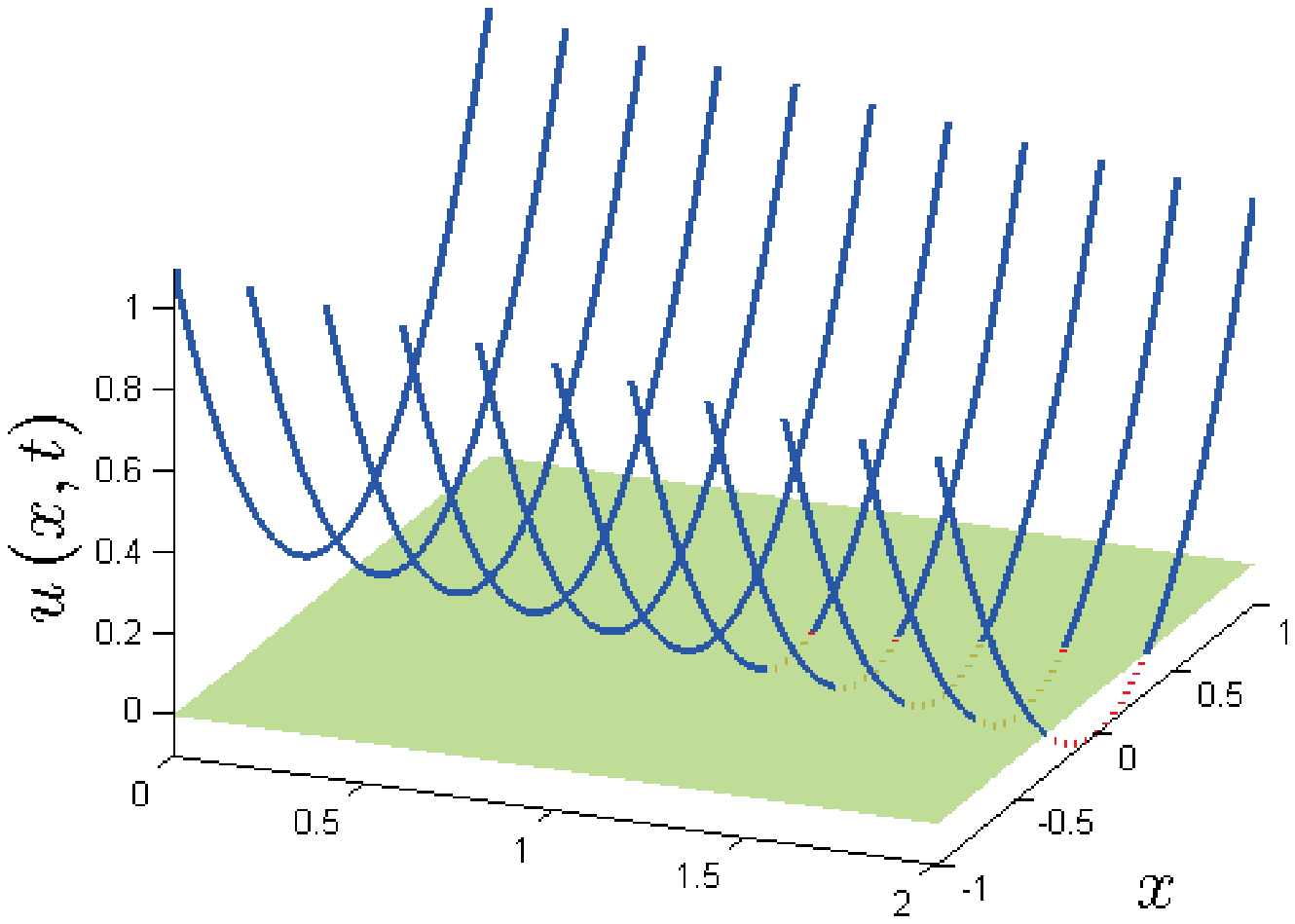}&
\includegraphics[scale=0.4]{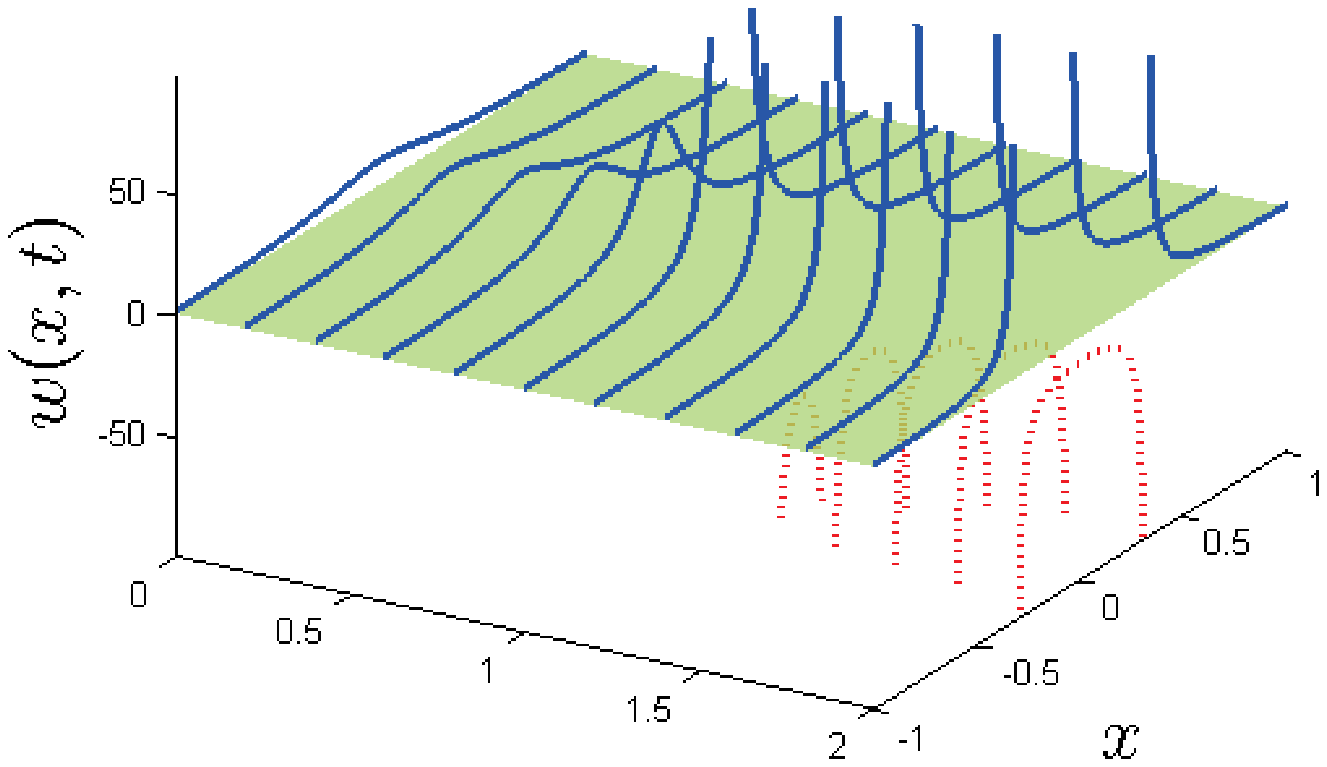}\\
(a)& (b) \\[8pt]
\includegraphics[scale=0.16]{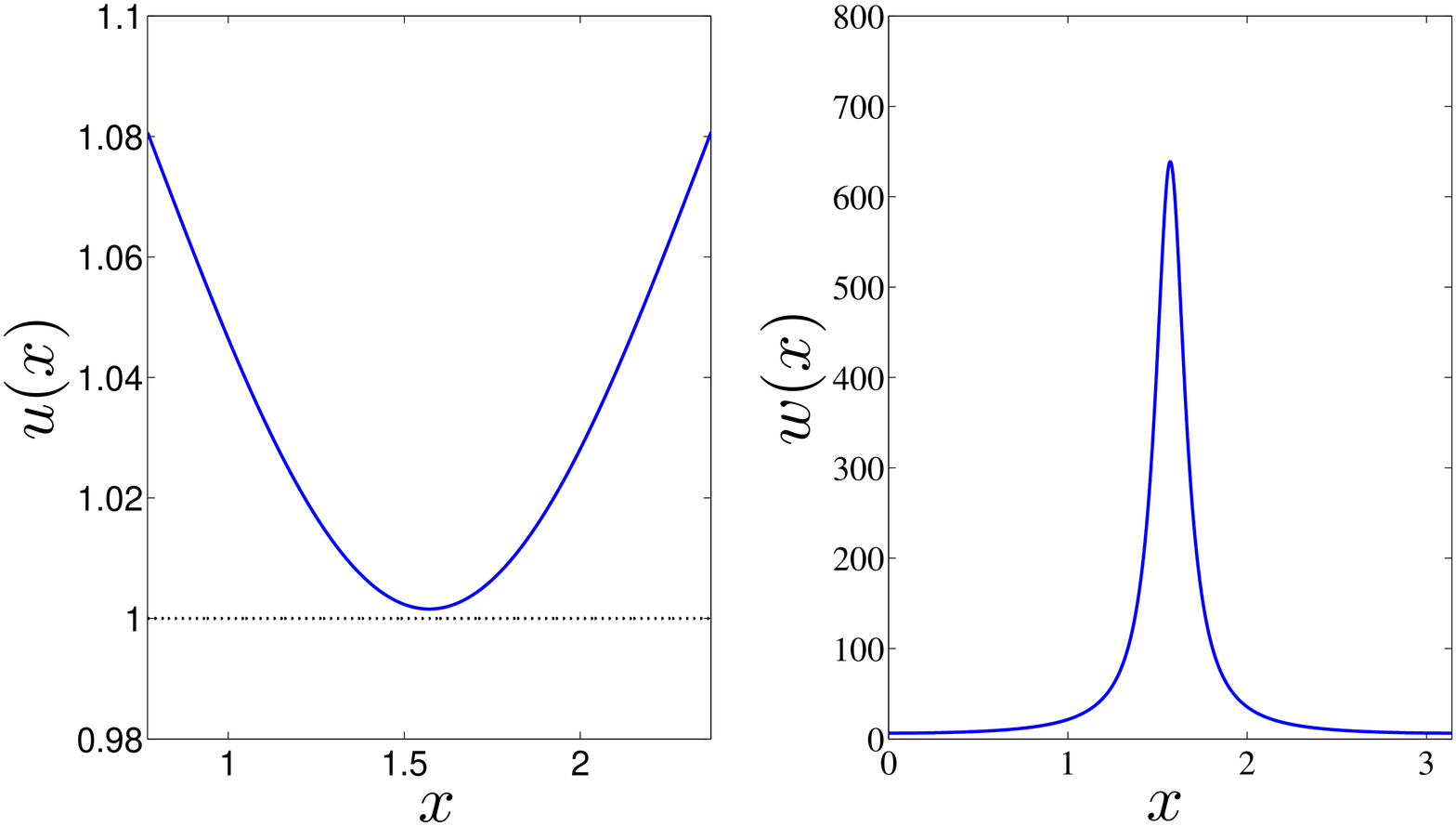}&
\includegraphics[scale=0.16]{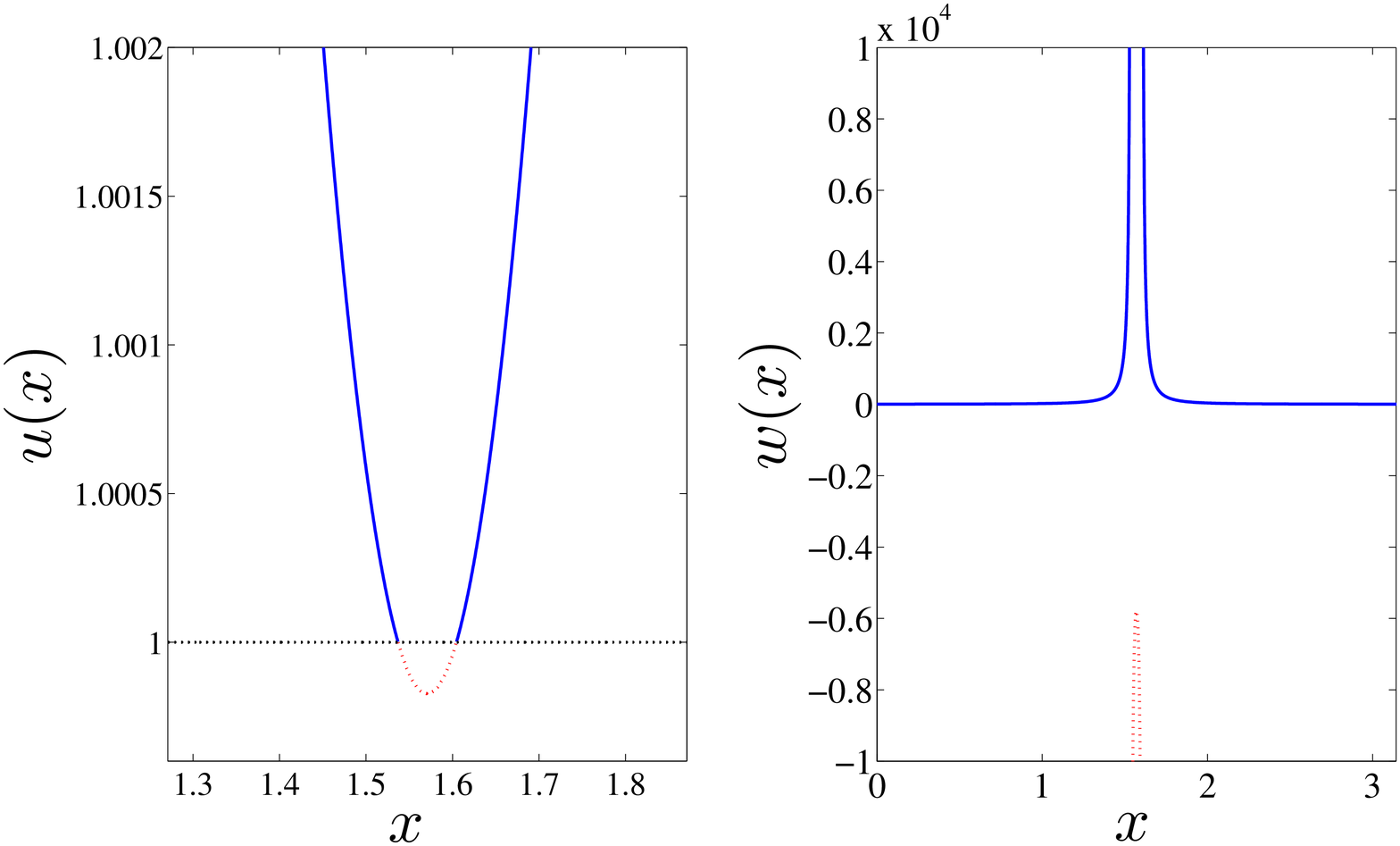}\\
(c)& (d) \\[8pt]
\end{tabular}
\caption{(Color online) A geometry for crossing infinity: (a,b) An idealized parabola   
$v(x,t)$ (a) shifting at constant speed without change of shape, crossing the $0$ level, (b) $w(x,t)=1/(v(x,t)-0)$.
(c,d) Solutions $u(x,t)$ of Eq.~(\ref{pdeqn8}) around their crossing of the $1$ level, (left panels), and corresponding
$w(x,t)=1/(u(x,t)-1)$ (solutions of Eq.~(\ref{pdeqn10})) (right panels).  (c): just before crossing, and (d): just after crossing infinity.}
\label{fig4}
%\end{center}
\end{figure*}

We now turn to a PDE example, illustrating one of the ways that
space dependence modifies the crossing of infinity.
%
%Scale equivariance now does not only involve the rescaling of
%an amplitude; it also involves the rescaling (typically dilation) of space,
%and self-similar collapse may lead to different post-collapse outcomes.
%
Motivated by $\frac{dx}{dt}= \pm x^2$, where $x^{-1}$ crosses infinity
at a single moment in time, we study the case where infinity is first reached 
in finite time, and then crossed continuously in time, but (in one spatial dimension) at
isolated points in space.
The geometry involved is illustrated in the top left panel of Fig.~\ref{fig4},
showing a graph of the function $v(x,t) = x^2 + (0.1-0.1t)$, a parabola
shifting its values downwards, at constant speed, but without change of shape.
Initially
it is everywhere positive; the tip reaches the zero level set at $t^*=1$ and then
crosses it. 
The function $w(x,t)=1/v(x,t)$ is shown in the top right
panel of Fig.~\ref{fig4}: it reaches
infinity at $t^*=1$ and subsequently crosses at two points that move apart
as dictated by the motion of the parabola.

We can then agree that the waveform ``returns from minus infinity''
between these two crossing points as 
%$w(x,t)=1/(u(x,t)-0)$ 
the definition of $w(x,t)$ formally suggests.
Extension to higher-dimensional geometries (e.g. a paraboloid initially touching a plane
at a point and then crossing it along a closed curve, such as a circle
that ``opens up'' in time starting at the initial point in two spatial dimensions) can also naturally be envisaged.
The key observation is that the evolution
of $w(x,t)$ actually involves a (potentially asymptotically) self-similar collapse near the crossing of
infinity.
%%
%% with the exponent of the self-similarity dictated by the speed
%% of the tip of the parabola.  HOW ARE THEY RELATED ?
%%
%% PGK: Consider your fictitious example x^2 + 0.1*(1-t). If 
%% you write w=1/((x^2+0.1 (t^*-t))=(1/(t^*-t)) 1/(\xi^2+0.1).
%% where \xi=x/\sqrt{t^*-t}
%% Then it's clear that the 0.1 does not dictate the exponent 
%% of the self-similarity...
%%
%%YGK
%% does the acceleration of the parabola tip 
%% have something to do with the exponent changing ?
%%
%% PGK: I am not sure what exponent changing we are 
%% talking about now. At the "center" x=0, the tip is
%% growing as $1/(0.1 (t^*-t))$ so it transparently gives
%% rise to the exponent.
This suggests that, upon detection -on the fly- of such an asymptotically self-similar collapse 
{\em and estimation of the associated exponents} (see below) 
for a ``bad'' $w(x,t)$ PDE, a search
for a ``good observable'' $v(x,t)$ be performed, so that a
conceptual and computational program analogous to that of 
the previous Section on ODEs may be carried through
to obtain, and work with, a "good PDE" in the vicinity
of the collapse point.

The simple, linear 
%reaction-diffusion 
equation
%% PGK: Do we REALLY want to write it is a REACTION
%% diffusion eqn. when there is no reaction -- just extinction ??
\begin{eqnarray}
u_t=u_{xx}- u
\label{pdeqn8}
\end{eqnarray}
provides an "engineered", yet transparent and analytically tractable illustration of the 
relevant ideas and hence will be used as our
workhorse in what follows.
Generic initial data in this
well-posed, linear model decay and concurrently spread,
asymptoting to $u(x,t)=0$  at long times.
We select an arbitrary level set $u^* = r > 0$ and a modified variable
$w \equiv \frac{1}{u(x,t)-r}$ to study level set crossings; for initial conditions
everywhere above $r$, and on its way to zero,
$u(x,t)$ will cross the level set $r$ so that $w(x,t)$ will cross {\em the level
set at infinity}.
%
%Having in mind the type of inverse
%transformation that we used in the case of ODEs,
%we wish for the transformation to be singular.
%To that effect we consider a modified variable $v=u-r$
%for which the PDE dynamics straightforwardly becomes:
%\begin{eqnarray}
%v_t = v_{xx} - v -r
%\label{pdeqn9}
%\end{eqnarray}

The ``bad'' PDE for $w(x,t)$ reads
\begin{eqnarray}
w_t = w_{xx} - \frac{2}{w} w_x^2 + w + r w^2.
\label{pdeqn10}
\end{eqnarray}
An auxiliary tool for the analysis of (asymptotic) self-similar collapse
in such equations is the so-called MN-dynamics~\cite{betelu,siettos};
%(see also details in the Appendix);
a dynamic renormalization scheme rescaling space, time and
the amplitude of the solution so that the self-similar solution becomes
a steady state in the ``co-exploding'' frame, i.e., the frame 
factoring out the symmetry/invariance associated with the 
(potentially asymptotic) self-similarity. 
This formulation is presented as a separate, detailed Section in the Appendix
for completeness.
In that Section, both the general case, and the special
example of Eq.~(\ref{pdeqn10}) are treated.
From this formulation we can infer that $w \sim 1/(t^{*}-t)$,
which, in turn, suggests the choice of a ``good variable''
used below.
%;
%in this frame a self-similar solution is sought as
%\begin{eqnarray}
%w(x,t) = A(\tau)f(y,\tau), y=\frac{x}{B(\tau)}.
%\label{selfsim1}
%\end{eqnarray}
%
%Putting \ref{selfsim1} into \ref{pdeqn10}, the MN dynamics read:
%\begin{eqnarray}
%f_t+G(\tau)(f+\frac{1}{2}y f(y)= f_{yy} - \frac{2}{f} f_{y}^2 + rf^2 + \frac{f}%{A(\tau)}
%\label{MN1}
%\end{eqnarray}
%where $\tau_t=A=\frac{1}{B^2}$, $\frac{B_{\tau}}{B}=-\frac{G}{2}$.
%
%Solving this equation with the algebraic constraint ....
%and with the boundary conditions .....
%eventually leads asymptotically to a steady state ...
%%%
%%% YGK
%%%  Itako, can I rely on you to asnwer this ?
%%%   if we do MN dynamics in an asymptotically self similar PDE
%%%    will it go to the self-sim, or do we have to go in "in person"
%%%     and find the offending terms and take them out ??
%%%
%%%%%%%%%%%%   so can we say that running the MN dynamics with the
%%%%%%%%%%%%   constraints will get the right stuff by itself at st. state ?
%%%
%%% PGK: No, because I simply didn't write this and
%%% I didn't write it this way in the Appendix. 
%%% I wrote it in a self-consistent general way, with this case
%%% as a transparent systematic case example without keeping
%%% in some places \tau and in some places t.
%%% Moreover, we have not discussed algebraic constraints for this
%%% problem, nor seriously worked on MN boundary conditions. You want me to make
%%% arbitrary statements that we have not tested and that are 
%%% irrelevant to the main story ? Why ???
%%%

In our illustrative example we use Neumann BC in $[0, \pi]$ and initial conditions
$u(x,0) = a \cos(2 x)+c$  (here $a =0.4, c =1.5 $,  
so that the solution $u(x,t)$ of  Eq. (\ref{pdeqn8}) reads:
\begin{eqnarray}
u = 0.4 \exp(-5*t) \cos(2*x)+1.5 \exp(-t),
\label{pdeqn9}
\end{eqnarray}
and we choose $r=1$.
We do not, however, pre-assume such knowledge of $u(x,t)$ since the equation we have to solve is the ``bad'' (focusing)
$w(x,t)$ equation, i.e., Eq.~(\ref{pdeqn10}); 
our MN framework applied to the focusing of the $w(x,t)$ evolution then suggests  that a good observable is $v(x,t) \equiv w(x,t)^{-1}$, 
a variable that will simply be crossing $0$ and thus the
"good PDE" would simply be
\begin{eqnarray}
v_t = v_{xx} - v -r
\label{pdeqn9a}
\end{eqnarray}
The bottom panels of Fig.~\ref{fig4} 
show representative instances just before and just after
the initial encounter of the $w(x,t)$ profile 
with infinity in both its "good" $v(x,t) \equiv (u(x,t)-1)$ and its "bad" 
$w(x,t)$ incarnations,
in the spirit of Figures 2a and 2b. 
The approach to infinity for
$w(x,t)$ is indeed {\em asymptotically self-similar}, as explained in the 
Appendix.
As we approach the event, an inverted bell-shaped profile comes close to, touches, and then
starts crossing through $r=1$ in the variable $u$, or crossing through $0$ in the variable
$v \equiv u-r$, or equivalently crossing through $\infty$ in the variable $w$.
%
%This is the computation we would like 

%Part of the reason for selecting this particular example is
%its full analytical tractability. Eq.~(\ref{pdeqn8}) possesses
%an explicit solution over the infinite line for arbitrary
%initial data in the form:
%\begin{eqnarray}
%u(x,t)=\frac{e^{-t}}{\sqrt{4 \pi t}} \int_{-\infty}^{\infty}
%e^{-\frac{(x-y)^2}{4 t}} u_0(y,0) dy
%\label{pdeqn10a}
%\end{eqnarray}
%for initial data given by $u_0(x,0)$. Then, $v$ can be obtained
%by a mere shift by $r$, while $w=1/(u-r)$.

%
Recall that our goal is to seamlessly
carry out the computation without our numerical code ever realizing that
(some part of) the solution is becoming indefinitely large.
To achieve this,
as the bad PDE solution grows towards infinity, it
is adaptively tested, with a user-defined threshold for 
(local, asymptotic) self-similarity, i.e., for growth 
according to a (potentially approximate) power law.
When this is numerically confirmed, a suitable power law transformation is devised with the 
numerically estimated similarity exponent; in the above example
the detected exponent is $-1$ and so the transformation is $v=w^{-1}$. 
The easiest way to realize the right observable in this case is to consider
uniform initial conditions - then the PDE reduces to an ODE that asymptotically
explodes as the $\dot{w}=w^2$, suggesting the $v=w^{-1}$ change of observables.

\begin{figure*}[h!]
\begin{center}
\includegraphics[scale=0.5]{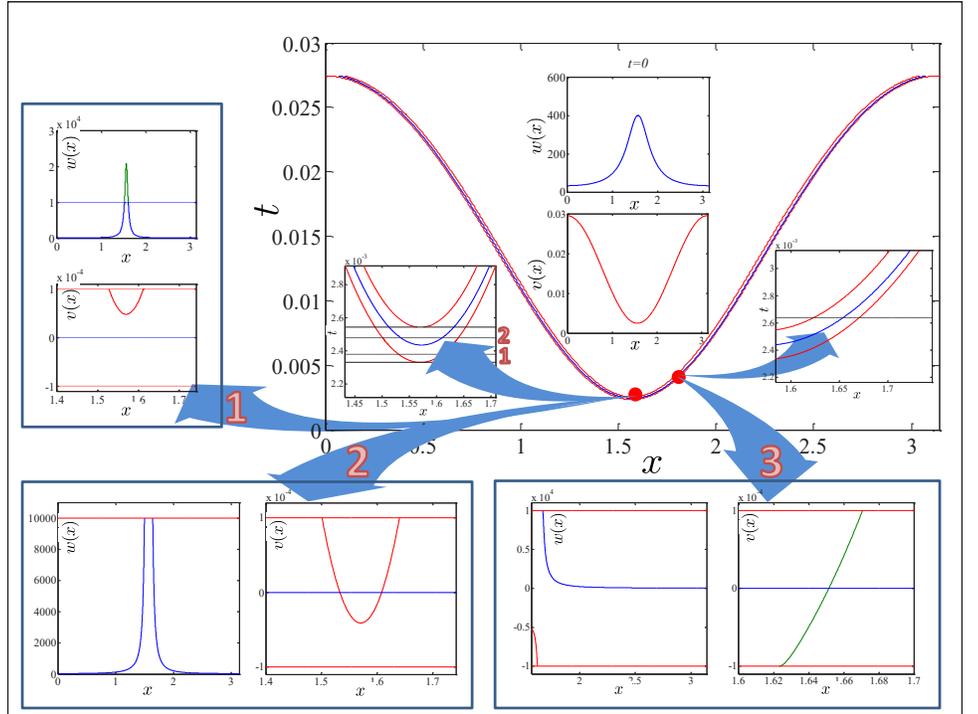}
\caption{Schematic of the computation. The main background figure (top right) shows the locus (in blue) of the points in space-time where $w$ becomes infinite. 
The two accompanying red curves
define the bands within which we solve the good equation for $v=1/w$.
The central top inset is the initial bad equation profile at $t=0$.  
The central bottom inset is the initial profile of our original motivating linear problem for $u(x,t)$.
The right inset shows a band crossing infinity in space-time, while the panels (3) show the solution of the good equation ({\em inside the band}, on the right) and of the bad equation ({\em outside the band}, on the left).
The left inset is more involved since it encompasses two computational ``eras''.
The first (1), takes place after we have crossed our selected high level (here $10^4$) for $w$ defining the buffer region boundary,
but have not yet crossed infinity; both the bad ($w$) and good ($v$) solutions are shown.
The second (2) takes place after $w$ crossed infinity, but has not (in absolute values) receded below the high level defining the buffer regions. Again both the bad ($w$) and the good ($v$) solutions are shown. 
The link to a computational movie of this evolution can be found at https://www.dropbox.com/s/1vag91zuwhjku1a/PDEPassInfinity.avi?dl=0} 
\label{fig5}
\end{center}
\end{figure*}

Importantly, the transformation
has to be performed -and the "good" solution sought- over an entire {\it spatial interval(s)} surrounding the approaching singular point(s). 
This suggests the following procedure, illustrated schematically in Fig. 3: 
\begin{itemize}
\item Upon detection of approach to infinity -as the ``tip'' of the collapsing waveform
grows beyond a sufficiently large value- at a given point or points inside the
computational domain, we split the domain in 3 regions:
(a),(c) regular ones to the left and to the right of the
growing tip, where the original ``bad'' equation for $w$ is being solved; and (b)
a new "singular"  one, in the middle, where instead of solving the equation
for $w$, now the equation for its singularly transformed variant, 
the "good" equation for $v=w^{-1}$ is solved instead. 
The latter
transformation is selected to comply with either the self-similar
analysis on the theoretical side, or the identified power law of amplitude
growth on the numerical side.
These equations are linked by continuity of the (transformed) observables
at the domain boundaries, and standard domain decomposition numerical
techniques are used~\cite{dupont}.
%%%
%%%  the reference is
%%%  C. N. Dawson, Q. Du and T. F. Dupont 
%%%  A finite difference domain decomposition algorithm for numerical solution
%%%  of the heat equation
%%%   TRF report 
%%%  and the good reference is 
%%%      Math. Comp. 57 (195) pp.63-71 (1991)
%%%
%
The good  equation simply crosses zero
rather than crossing infinity, as in the ODE case.

\item Once zero is crossed, the initial single crossing point 
in the family of case examples under consideration (in one
spatial dimension) "opens up" into {\em two}
infinity crossings (one can visualize two waves that propagate in opposite directions,
one to the left and one to the right) - two zero-level-set crossings for the "good"
equation.
These crossings are quantified, for our example, in Fig.~\ref{fig5}. 
They are bordered by the computed locations of a "high enough" absolute level
(here $10^4$) for asymptotic self-similarity. 

To deal with the two new crossings computationally over time, the central region is subsequently
split into 3 regions. The two outer ones are our "singular buffers",  surrounding, and in some sense "masking" the infinity crossings to the left and to the right.
But now they are separated by another, inner "regular" interval, where we can again
solve the original "bad" equation since in here it is again sufficiently
far from infinity.
%We recap the procedure: close to -but before- the collapse event the domain
%is split into three regions, two regular ones, where we solve the
%original bad equation, and an inner ``singular" one, where we
%solve the transformed ``good" equation, with continuity (modulo the
%transformation) at the transitions.
%Post-collapse we have two options (the reason for the second one will be discussed
%further below). 
%
%The simplest is to keep three regions, and keep solving the
%"bad" equation in the two outer ones, and the "good" in the middle.
%
%Alternatively, it is also possible, and as we will argue, even desirable,
Thus, post-collapse, we partition the domain into five regions, three regular ones -the two outer ones, 
and the innermost, for the
"bad equation"- and then two "singular buffers" for the transformed "good" equation,
one around the left zero crossing and one around the right zero crossing of $v$,  that correspond to the two infinity crossings of the bad equation for $w$.
\end{itemize}

A nontrivial aspect of the computation is the "gluing" between
the regular regions and the singular buffers. 
Our numerical scheme here is a simple one following~\cite{dupont}: for a finite difference discretization in space, (a) an explicit forward Euler time step is performed at the interface points, which provides the interior boundary conditions for the next time step (the next ``computational era"), while (b) an implicit Euler time step is adopted to solve the "good" and "bad" PDE within the three (or the five) domains. The scheme can be modified to allow for different space and time steps in the different domains till the next computational era, when the new interface points will be detected, the new "interior" BCs will be computed, and the new set of BVPs resulting from the implicit timestepping in each domain will be solved.
To recap the essence of the algorithm, by solving the "good" equation
inside the buffer regions (and following the motion of the buffers on the fly), we ensure
that the numerical simulation is never plagued by the indeterminacy
associated with approaching/touching/crossing infinity.

We now explore two more notions that were also examined in the ODE context. 
The first is compactification: self-similar PDE dynamics, 
although crossing through infinity,
can simply be compactified as evolving over a circle or -better- on a sphere.
We perform this step for a general 1D PDE (self-similar) solution of
the form --assuming independence from $\tau$ (see the MN formulation
in the Appendix)--:
\begin{eqnarray}
u(x,t)= \frac{1}{(t^{\star}-t)^r} f(\xi),
\label{pdeqn12}
\end{eqnarray}
where $\xi$ is a self-similar variable, e.g. $\xi=x/(t^{\star}-t)^q$.
If the $\max(f)>1$, we can define $\tilde{f}=f/\max(f)$ and
subsequently drop the tilde in Eq.~(\ref{pdeqn12}) ensuring
that $|f| \leq 1$.  
We can then rewrite Eq.~(\ref{pdeqn12}) as:
\begin{eqnarray}
(t^{\star}-t)^r u = f \Rightarrow
\left((t^{\star}-t)^r + u\right)^2 -  \left((t^{\star}-t)^r - u\right)^2 =
4 f(\xi).
\label{pdeqn13}
\end{eqnarray}
Then, upon suitable definition of the variables,
we can have $X=\left((t^{\star}-t)^r - u\right)/\left((t^{\star}-t)^r + u\right)$
and $Y=2 \sqrt{|f(\xi)}/\left((t^{\star}-t)^r + u\right)$, in which
case $X^2+Y^2=1$. 
In these variables, at every moment in time
the trajectory can be thought of as compactified along a circle.
However, as the circle itself represents an invariant shape,
in this representation we cannot straightforwardly visualize the trajectory's dynamics;
for this reason, we next compactify the dynamics {on a sphere}.
We define $g^2=1-f^2$ and we can then write using the above
variables $(g X)^2 + (g Y)^2 = g^2=1 - f^2$, which can be reshuffled
to read:
\begin{eqnarray}
\left( g  \frac{(t^{\star}-t)^r - u}{(t^{\star}-t)^r + u} \right)^2
+ \left( g  \frac{2 \sqrt{f}}{(t^{\star}-t)^r + u} \right)^2
+ f^2 =1.
\label{pdeqn14}
\end{eqnarray}
Choosing the three terms of the left hand side of Eq.~(\ref{pdeqn14})
as $(X',Y',Z')=(g  \frac{(t^{\star}-t)^r - u}{(t^{\star}-t)^r + u},
g  \frac{2 \sqrt{f}}{(t^{\star}-t)^r + u},f)$, we observe
that the dynamics can be seen as evolving along the surface of
a sphere. 
This compactification once again underscores the
possibility to think of infinity as a regular circle
(rather than point, as is the case for ODEs), a level set that is
crossed by the PDE solution evolving along the surface of the sphere.

\begin{figure*}[tb]
\begin{tabular}{cc}
\includegraphics[scale=0.4]{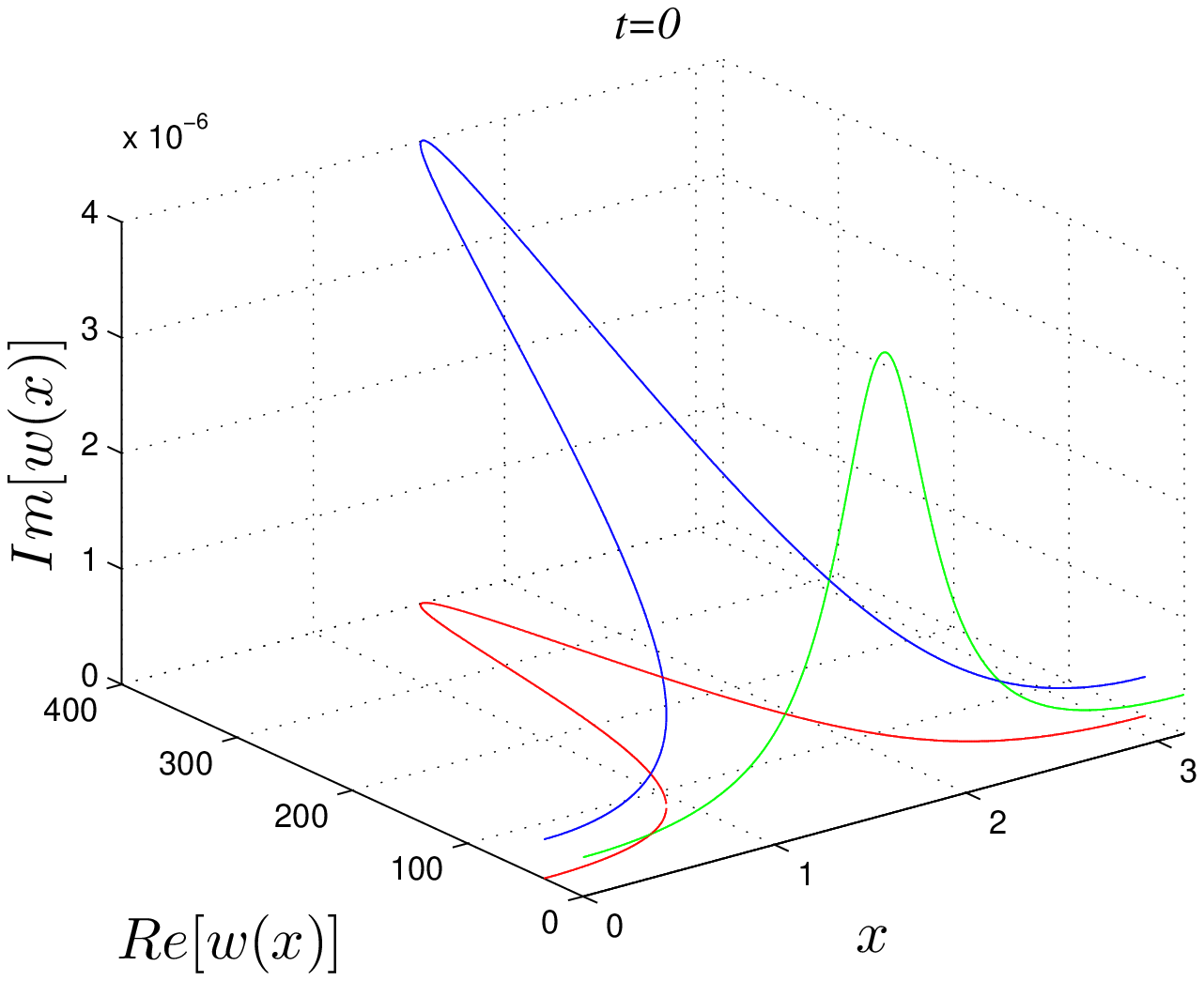} & \includegraphics[scale=0.4]{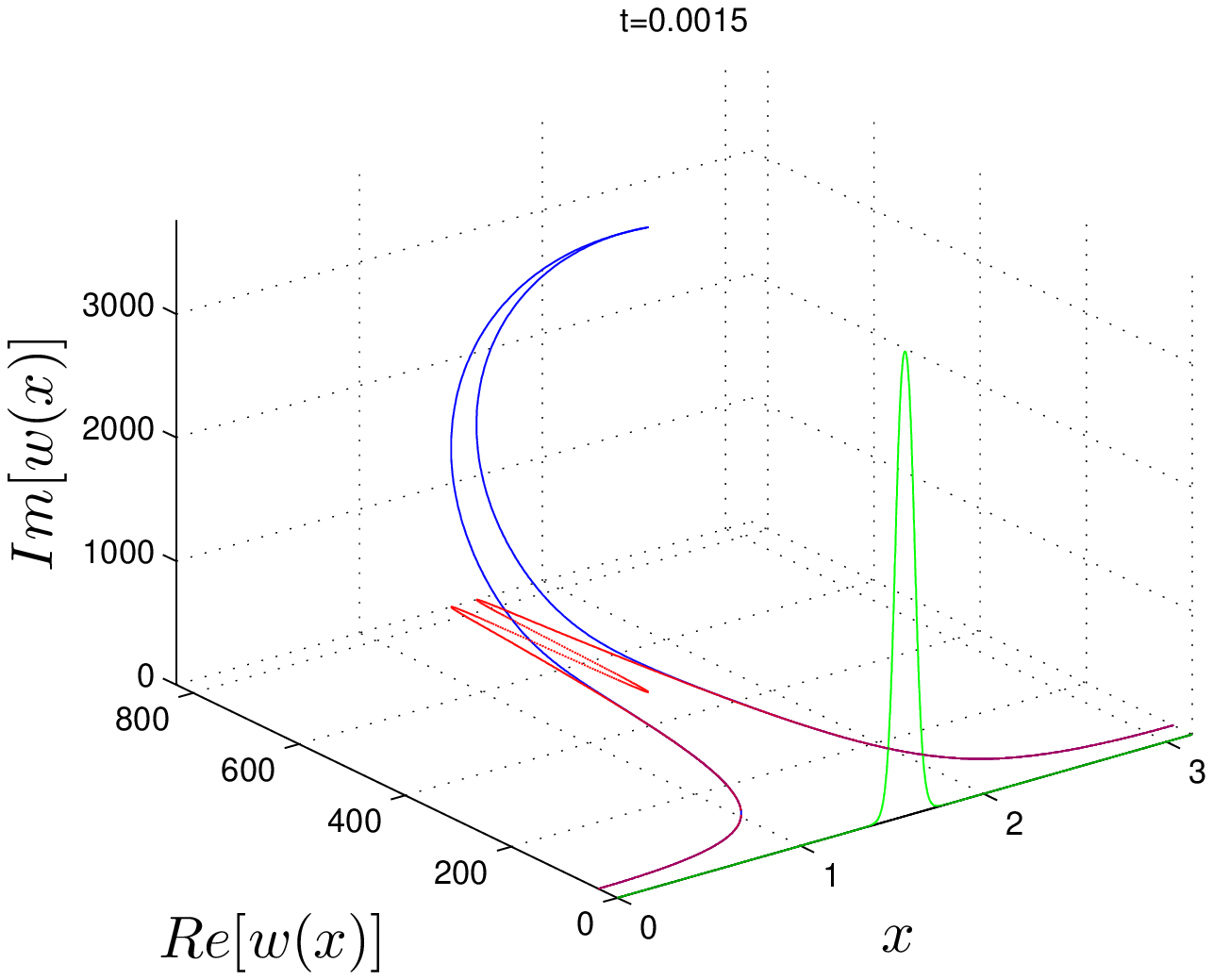}\\
(a)& (b) \\[8pt]
\includegraphics[scale=0.4]{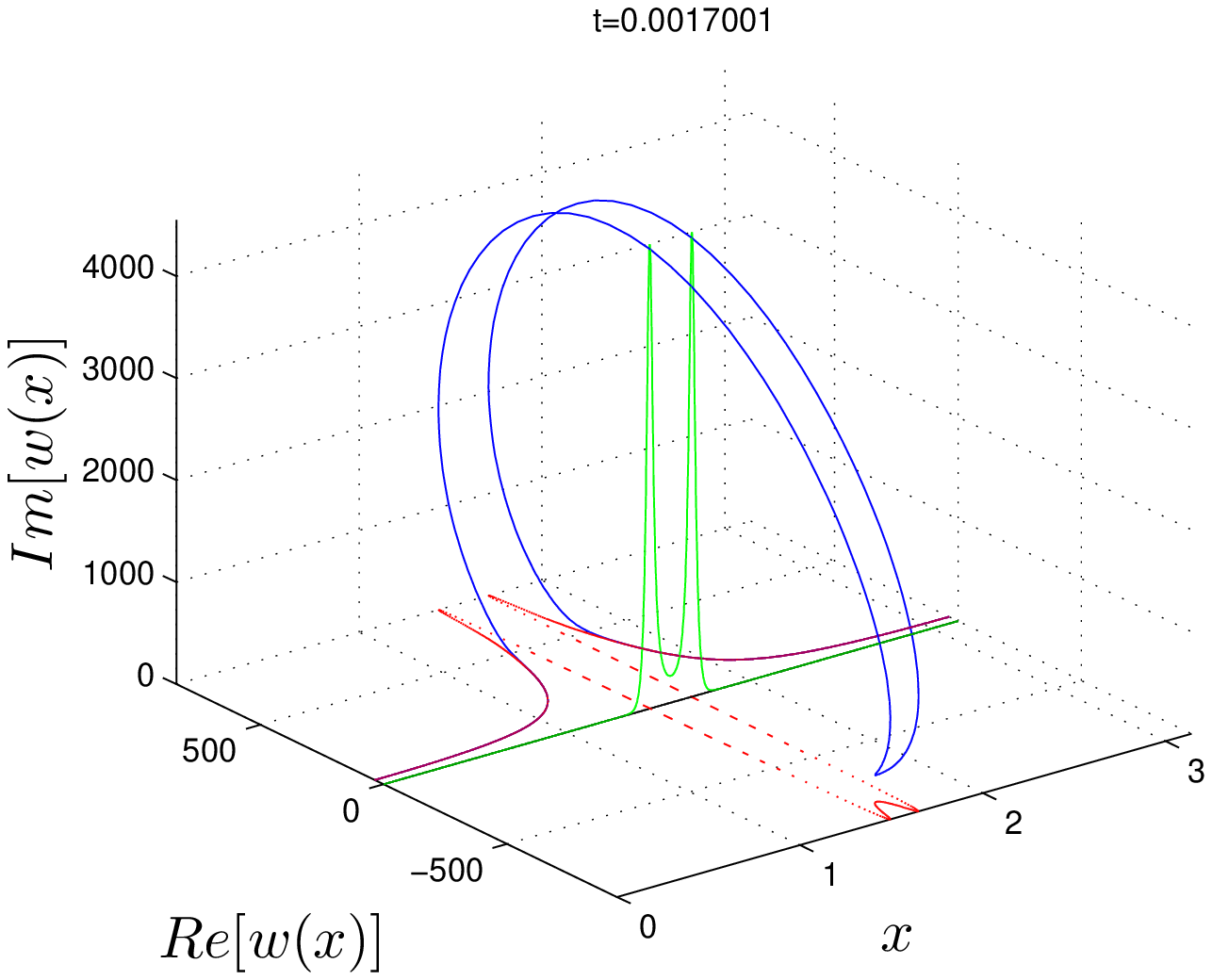} &
\includegraphics[scale=0.4]{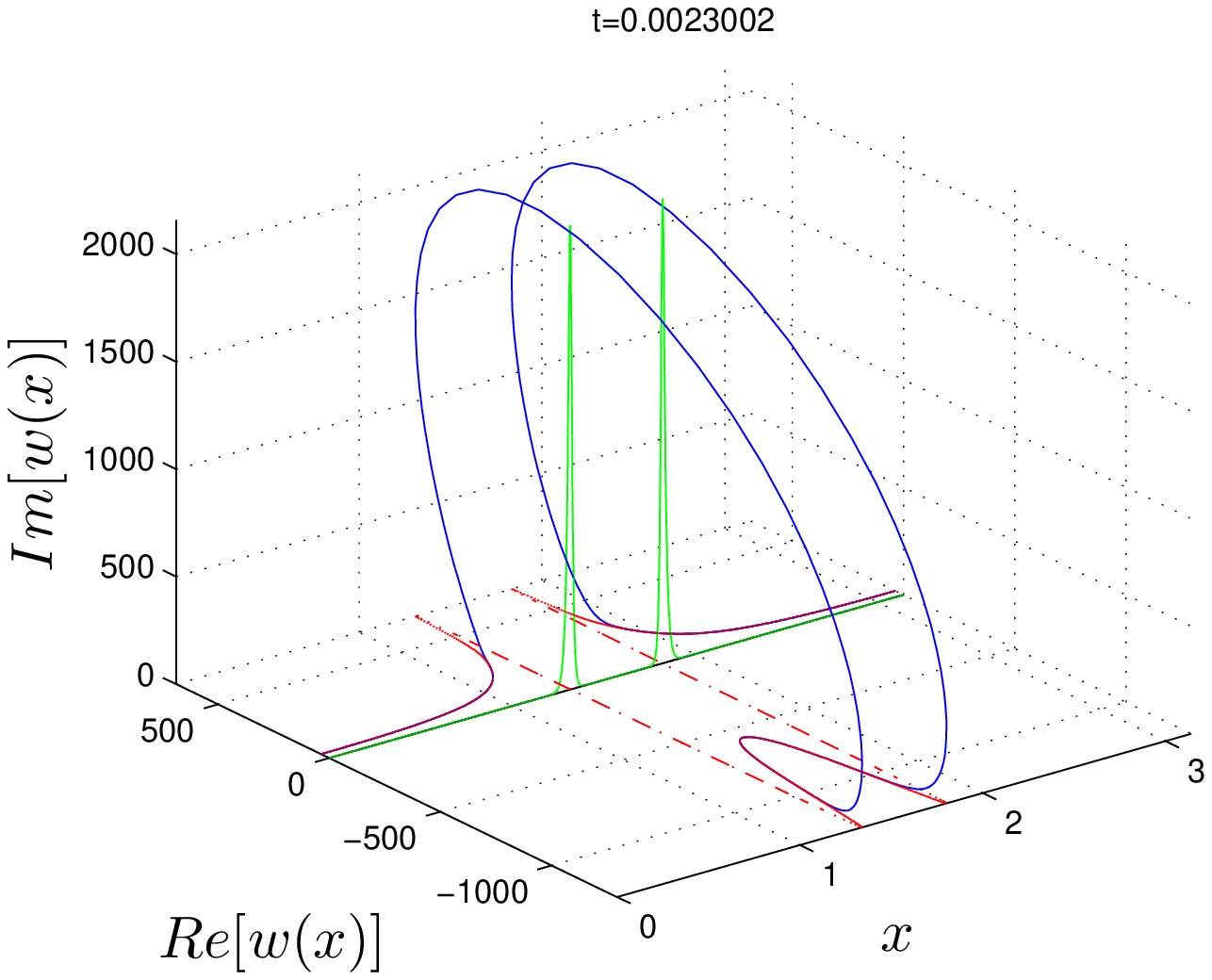}\\
(c) & (d) \\[8pt]
\end{tabular}
\caption{(Color online) The figure shows the solution of the complex equation associated with the real and imaginary parts of Eqs.~(\ref{pdeqn15})-(\ref{pdeqn16}) at various time instances (a) $t=0$, (b) just before the collapse time, where the complex solution starts ``going back on itself'' rather than collapse (as it would for purely real initial data), (c), (d) Snapshots of the solution ``reemerging" on the other side (the negative real axis end of the complex plane). Notice the visual similarity with the images in Fig. 2.
The link to a computational movie of this evolution can be found at https://www.dropbox.com/s/rxekr7umxa5b3ko/complexdynamics.avi?dl=0}
%The left panel is depicting the complex solution profile in
%terms of the spatial distribution of its real and imaginary
%part at times $t=0$, $t=1$, $t=1.5$, $t=1.8$, $t=2$. Here, we see the tendency of the real part towards collapse. The middle panel one is for times
%$t=2$, $t=2.1$, $t=2.2$, $t=2.3$; in each panel we also show the end profile of the previous panel for comparison purposes. Here, we see the complexification leading to the regularization of the solution by leading it to ``go back on itself'' rather than to truly collapse (as it would for purely real initial data). Finally, the third panel
%is for $t=2.3$, $t=2.4$, $2.45$, and $t=2.6$ where the solution eventually reemerges on the other side (the negative real axis end of the complex plane).}
\label{fig6}
%\end{center}
\end{figure*}

Finally, as in the ODE case, we discuss the possibility of complexifying the
model in order to understand, as a limiting case, how  infinity is crossed for
purely real initial data, while it may be avoided
(regularized) upon initialization with complex
initial data. 
Using the complex decomposition for
$w(x)=a(x)+ i b(x)$ in Eq.~(\ref{pdeqn10}), one can obtain the pair of real and
rather elaborate looking equations:
\begin{eqnarray}
a_t &=& a_{xx} -2 a \frac{a_x^2-b_x^2}{a^2+b^2} - 4 b \frac{a_x b_x}{a^2+b^2}
+a + r (a^2-b^2)
\label{pdeqn15}
\\
b_t &=& b_{xx} + 2 b \frac{a_x^2-b_x^2}{a^2+b^2} -4 a \frac{a_x b_x}{a^2+b^2}
+ b + 2 r a b.
\label{pdeqn16}
\end{eqnarray}
We expect that the presence an imaginary part in the initial
data may avoid collapse in analogy with Figure ~\ref{fig2}.
Given the quadratic nature of the nonlinearity,
the quadratic ODE example is especially relevant; we expect
here to observe something similar but in a PDE form, having space
as an additional variable, over which the profile is
distributed (around the crossing ``tip''). 
Again, the tractability of our
example allows us, via the solution of Eq.~(\ref{pdeqn9a}),
to perform the relevant calculation analytically since at the level
of the equation for $v$ the complex model can be fully solved.
Then, assuming $v=c+ i d$, the variable $w= a + i b=1/v=1/(c+ i d)$
leads to $a=c/(c^2+d^2)$ and $b=-d/(c^2+d^2)$, and obtaining
$(c,d)$ explicitly, the same can be done for $(a,b)$.
This program is carried out in Fig.~\ref{fig6}; see also the
relevant movie at:
https://www.dropbox.com/s/rxekr7umxa5b3ko/complexdynamics.avi?dl=0. 
We reconstruct
analytically the spatial profile of the real and imaginary
parts of $w$ at different moments in time provided in the caption.
While the profile tends towards collapse in the real part of
the variable (and would go all the way to collapse and the
dynamics of a Figure like~\ref{fig4} for purely real initial data),
the imaginary part, in analogy to the dynamics of Fig.~\ref{fig2},
but in a distributed sense around the tip, eventually takes over.
As it does so, it forces the solution filament to ``turn around
on itself'' in a spatially distributed generalization of
the ODE of Fig.~\ref{fig2}. 
Finally, the solution appears to re-emerge from
the other side, practically extinguishing its imaginary part,
and having avoided the crossing of infinity. 
This illustrates how complexification, even in the case of the PDE, results
in the avoidance of collapse and the regularization of the model;
the collapsing real case is a special limit case of the more
general complex one.
%

%% PGK: Here, as well as in the figure 
%% we need to give a link to the relevant movie...
%% While the story in this case is well phrased by
%% the snapshots the movie would do even better
%% In the case of the previous Fig. (less clear and
%% more busy) the movie will even go a longer way...

\section*{Concluding remarks and future challenges}

%%%YGK_NOW  you need to talk with me here. There IS a bigger picture
%%%  and we are ignoring it.  There is a chance we will appear blind
%%%  for having missed it.
%%%

%% PGK: I don't know what the bigger picture is in your mind.
%% I can think a ton of extensions in all sorts of 
%% directions (and I have already 
%% written down nearly half a dozen that are stoically 
%% awaiting their turn after this finishes and is broken out
%% (NLS defocusing, NLT-PT-Ziad, the
%% 2d collapse on a ring, the connection with the algorithm
%% problem and Davydenko, Burgers-and-complex-Burgers, etc.) 
%% while nothing happens... For me, these are part of a bigger
%% picture and program which I somehow do hope to pursue.
%% From there on, needless to say, you are welcome to 
%% write anything that you see fit (and if you 'd like me
%% to proof read it, I 'll try to do so).
We attempted to address here a few prototypical cases
of a spectrum of problems arising in both ordinary
and partial differential equations, so as to deal with the emergence
of infinities during the 
evolution of the relevant models.
In a number of cases the model
at hand will become physically inaccurate, and will need to be suitably
modified as these singular points are approached; if not,
our questions may be relevant for the physical realm. 
In any event, the questions are of particular relevance towards the
mathematical analysis and numerical computation of the models
at hand. In that light, we argued that it is possible in our context to perform
singular transformations on demand, that may sidestep
-through the help of a "good equation"- the computational difficulties 
associated with infinities,
rendering them tantamount to the crossing of a regular
point such as zero.
For ordinary differential equations, once the crossing has transpired, 
one can safely return to the 
original "bad" equation and continue the dynamics from there 
(until possibly a new infinity is approached).

In the case of partial differential equations, the scenario
at hand is more complex. There, the solution is distributed in space,
and hence we assume and have analyzed the setting where a
(generically assumed to be parabolic; see the relevant discussion in the 
Appendix) tip of a waveform approaches infinity. 
We have discussed in detail a scenario of initially touching infinity
and then crossing it. Suitable computational "buffers" need then to be devised, where the 
detected singular transformation allows us to locally
re-interpret (for computational purposes) the crossing of
infinity as the crossing (in a transformed space) of a regular point, such as zero.
These buffers need to be in constant and consistent
communication, through appropriate continuity conditions, with the rest of the
computational domain (the "rest of the world").
Typically, the buffers are defined by the location at which the
solution takes on a sufficiently large (absolute) value - say $10^4$ to the left and right
of the growing tip in the pre-crossing regime,
or, say, $10^4$ till $-10^{4}$ on the left, and $-10^4$ till $10^{4}$ on the right
in the post-crossing regime.
Our computational findings were
complemented by a compactification approach, supporting the argument
that infinity can be addressed in the same way as a regular point or
a regular level set along the orbit.
At the same time, a complexification of the model was observed to provide
a regularization of the original real dynamics, avoiding the collapse
of the latter and offering insight on how collapsing orbits can be 
envisioned as limiting scenarios of nonlinear dynamical 
systems within the complex plane.

Naturally, there are numerous directions of interest for potential
future studies. Clearly, exploring additional examples and
examining whether the ideas can be equally successfully applied
to them is of particular relevance. In the context of ODEs,
this is especially relevant as regards vector/multi-dimensional
systems.
A related, especially important part in the realm of ODEs
is that of convergence of algorithms e.g. to fixed points
(or extremizers) of functions. Recall that in such cases,
a concern always is whether the code may diverge along the way, rather
than reach a root (or an extremum). Our approach  can be used
to devise algorithms with the ability to systematically bypass
infinities during the algorithmic iterations- 
and such a "boosted" algorithm may be useful 
towards achieving enhanced, possibly global convergence to the roots 
(or extrema) of a function. This is particularly interesting now
that continuous time versions of time-honored discrete algorithms
like Newton or Nesterov iteration schemes have become a research focus; see the related discussions in~\cite{deflation,candes,Jordan}.

On the PDE side, we are envisioning (and currently starting to explore)
a multitude of emerging aspects. For instance,
when a distributed waveform reaches infinity at a single point in space-time, different 
post-collapse outcomes are possible.
For example, an alternative
possibility to the infinity-crossing presented here has been argued to be 
that the solution may "depart" from infinity 
without crossing (the transient blowup 
in~\cite{galakt}), as in the case of the
%as in
%the case of the diffusion equation $v_t = v_{xx}$
%(we will now use subscripts to indicate partial derivatives).
%Another model that has a qualitatively similar feature is the
standard collapsing NLS equation discussed extensively
in textbooks~\cite{sulem,fibich}.
There, the crossing through infinity is precluded by the existence of conservation laws.
Past the initial point, it is
argued in~\cite{fibich,gaeta} that the solution
will return from infinity incurring a ``loss of phase''.
At the bifurcation level, the work of~\cite{siettos} offers a suggestion of how
the return from infinity manifests itself: there, a solution
with a positive growth rate was identified, that was
dynamically approached during the collapse stage.
Yet a partial ``mirror image" of that, with negative growth rate,
which presumably is followed past the collapse point
in order to return from infinity was also identified; see, in particular,
Fig.~1 and especially Fig.~2 of~\cite{siettos}.
%
%%%

It is also possible that such a "touch and return"
from infinity may occur without the loss of phase
as, e.g.,  in the recent work of~\cite{mussli}.
In examining a nonlocal variant of NLS
(motivated by ${\mathcal PT}$-symmetric considerations,
i.e., systems invariant under the action of parity and time-reversal),
Ref.~\cite{mussli} identified a solution that goes to
infinity in finite time that can be theoretically calculated;
subsequently this solution returns from infinity and
then revisits infinity again, in a periodic way, always
solely touching it and never crossing. This solution
is analytically available in Eq.~(22) of~\cite{mussli} and
the collapse times are given by Eq.~(23) therein; perhaps even
more remarkably, the model itself is integrable.
In this case, infinity is reached, subsequently returned from and
then periodically revisited.
%
%%%YANNISNOW
Such an observation would arise in our context if the
"original" PDE for $u$  (i.e., a variant of Eq.~(\ref{pdeqn8})) had a spatiotemporal limit cycle that
attained somewhere in space an extremal value $r$.
Then, $w(x,t) \equiv \frac{1}{u(x,t) - r}$ would feature
the above phenomenology.
%
%We would simply have to  "strategically" choose our arbitrary level
%set $r$ in the definition of $w(x,t) \equiv \frac{1}{u(x,t) - r}$  to be %the maximum local value.
%This example is perhaps the most intriguing to
%explore with the technology developed herein, and is
%deferred to an immediate future calculation.
Such cases where infinity is reached but not crossed
merit separate examination.
The same is true for solutions exhibiting entire intervals at infinity,
whose support progressively grows (or anyway remains
finite), bordered by moving "fronts";
here one may envision that the "good" equation develops
compacton-like solutions~\cite{rosenau}. 
%%%
%%% Compactons: Solitons with finite wavelength
%%%  Philip Rosenau and James M. Hyman
%%%   Phys. Rev. Lett. 70, 564 – Published 1 February 1993
%%%   

A related issue that may be worth exploring with such techniques
is the possibility of bursting mechanisms  (e.g. ~\cite{knob1,knob2}
involving heteroclinic connections with entire invariant planes at infinity)
and the associated emergence of extreme events in nonlinear
PDEs.
Generalizations of the techniques
developed herein to settings where, rather than $u(x,t)$,
$u_x(x,t) \rightarrow \infty$ (or this happens for other quantities associated with the dependent variable), as is, e.g., the case during the formation of shocks, should also be interesting to explore.
%It will be interesting to explore analogies to heteroclinic
%excursions to and from an ``invariant plane at infinity'' \cite{...,.....},
%a mechanism potentially underlying rogue wave formation in PDEs.
Effectively, our considerations here can be thought of as identifying and numerically
evolving the infinity level set  of the solution.  
Thus, a related interesting direction for future work could be to try to connect
the considerations herein with ones of level set methods~\cite{sethian1,sethian2},
adapting the latter towards capturing, e.g., the regions of the singular buffers.

%We are confident that 
%the mathematical and computational technology developed herein will be of
%value in elucidating such examples in the near future.

Equally relevant
are explicit examples similar to the one herein 
where multiple collapses may occur and
propagate. An intriguing such case is the {\it defocusing} scenario
of the nonlinear Schr{\"o}dinger equation,
\begin{eqnarray}
i u_t= u_{xx} -2 |u|^2 u
\label{concleq1}
\end{eqnarray}
which, in fact, has been
shown in~\cite{clarkson} to possess solutions such as
$u(x,t)=1/x$, or
\begin{eqnarray}
u(x,t)=\frac{2 x (x^2 + 6 i t)}{x^4 - 12 t^2}
\label{concleq2}
\end{eqnarray}
with propagating singularities at $x= \pm 12^{1/4} t^{1/2}$,
and
\begin{eqnarray}
u(x,t)=\frac{3 (x^8 + 16 it x^6 -120 t^2 x^4 + 720 t^4)}{x (x^8-
72 t^2 x^4 - 2160 t^4)}.
\label{concleq3}
\end{eqnarray}
It is obvious that to follow such dynamical examples, a methodology
bearing features such as the ones discussed above is needed in order
to bypass the continuously propagating singular points.
In turn, generalizing such notions to higher dimensions 
(e.g., a two-dimensional variant of the analytically tractable
example herein) and addressing
collapsing waveforms both at points, as well as in more complex
geometric examples such as curves~\cite{gavish} is of particular
interest for future studies.
It is tempting to explore 
whether the tools developed here may have something
to add in the way we analyze collapse in well-established PDEs like the 
Navier-Stokes, or even singularities arising in a cosmological context. 
%%%
%%% 
%%%
Several of these topics are under active
consideration and we hope we will be able to report on them in future publications.
%\clearpage

\section*{Appendix}
{\bf Parabola Self-Similar Crossing}

In the 1d case, starting from the assumption of a self-similar
solution approaching infinity, we can prescribe a generic ``unimodal''
profile of the form:
\begin{eqnarray}
u \sim \frac{1}{(t^{\star}-t)^a} f\left(\frac{x-x_0}{(t^{\star}-t)^b}\right),
\label{ceq1}
\end{eqnarray}
where $t^{\star}$ is the collapse time and $x_0$ is the point
around which the blowup solution is centered.

Then, around $x=x_0$ (and for $t \neq t^{\star}$), we can use a Taylor
expansion locally in the form:
\begin{eqnarray}
u \sim \frac{1}{(t^{\star}-t)^a} \left[ f(0) + f'(0) 
\frac{x-x_0}{(t^{\star}-t)^b} + \frac{f''(0)}{2} 
\frac{(x-x_0)^2}{(t^{\star}-t)^{2 b}} \right]
\label{ceq2}
\end{eqnarray}

Combining the powers and bringing the dominant power to the left,
we obtain that the field $v$, defined as
\begin{eqnarray}
v \equiv u (t^{\star}-t)^{a + 2b} \sim \left[ f(0) (t^{\star}-t)^{2 b} +
\frac{f''(0)}{2} (x-x_0)^2 \right],
\label{ceq3}
\end{eqnarray}
behaves like a ``regular'' field which crosses $v=0$ at $x=x_0$,
when $t=t^{\star}$. So, its dynamics should be that of a 
``rising parabola'', cutting through 0 at the critical time.
In Eq. (\ref{ceq3}), we also used the fact that $x=x_0$ was an
extremum (having in mind in particular a maximum) of the profile 
of the solution (hence $f'(0)=0$).

%\smallskip
{\bf MN-Dynamics}
%\section{MN-Dynamics}

As an auxiliary tool in our analysis, we will outline here
and utilize the so-called MN-dynamics~\cite{betelu}, i.e.,
the self-similar dynamical evolution of a PDE which is
collapsing towards a dynamical formation of a singularity.
This approach has been used in porous medium type equations,
as well as in dispersive (and conservative) NLS equations~\cite{siettos}
and is broadly applicable
to problems with self-similar growth (or decay).
To illustrate it in a general form, we consider
an evolutionary PDE of the form:
\begin{eqnarray}
u_t = {\cal L}[\partial_{\xi}] u + {\cal N}[u],
\label{pdeqn3}
\end{eqnarray}
By ${\cal L}$ here we designate the operator involving
derivatives (not necessarily a linear operator -- see also
the example below), while by ${\bf N}$ we designate the
local nonlinearity bearing operator.

Using the ansatz
\begin{eqnarray}
u=A(\tau) f(\xi,\tau); \quad \xi=\frac{x}{L(\tau)}, \quad \tau=\tau(t)
\label{mndynamics}
\end{eqnarray}
we introduce a new scaled system of coordinates, intended to be
suitably adjusted to the self-similar variation of the PDE solution.
$\xi$ is a rescaled spatial variable (taking into consideration
the shrinkage --or growth-- of the width), while $\tau$ is a rescaled
time variable, not a priori tuned, but which will be adjusted
so that in this ``co-exploding'' frame, we factor out the self-similarity,
in the same way in which when going to the co-traveling frame, we factor
out translation. This way, the self-similar solution resulting in this dynamical
frame will appear to be steady. Direct substitution of Eq.~(\ref{mndynamics})
inside of Eq.~(\ref{pdeqn3}) yields:
\begin{eqnarray}
\left[A_{\tau} f + A f_{\tau} -A \xi f_{\xi} \frac{L_{\tau}}{L} \right]
\tau_t = {\cal L}[\partial_{\xi}] f \frac{A}{L^a} + A^s {\cal N}[f]
\label{pdeqn4}
\end{eqnarray}
where $a$ and $s$ are powers tailored to the particular problem
(linear and nonlinear operators) of interest.
In order to match the scalings of the two terms of the right hand
side of Eq.~(\ref{pdeqn4}), as is required for self-similarity,
we demand that:
\begin{eqnarray}
\frac{1}{L^a} = A^{s-1} \Rightarrow A \sim L^{-\frac{a}{s-1}}
\Rightarrow G \equiv \frac{A_{\tau}}{A}= - \frac{a}{s-1} \frac{L_{\tau}}{L}
\label{pdeqn5}
\end{eqnarray}
Thus, the model can now be rewritten as:
\begin{eqnarray}
\left[G \left( f + \frac{s-1}{a} \xi f_{\xi} \right)
+ f_{\tau} \right] \tau_t= A^{s-1} \left({\cal L}[\partial_{\xi}] f
+ {\cal N}[f] \right).
\label{pdeqn6}
\end{eqnarray}
Demanding then that the time transformation be such that there is
evolution towards a stationary state in this co-exploding frame,
we remove any explicit time dependence by necessitating that
$\tau_t = A^{s-1} \sim L^{-a}$. Then, the stationary state in this frame
will satisfy:
\begin{eqnarray}
G \left( f + \frac{s-1}{a} \xi f_{\xi} \right)=
{\cal L}[\partial_{\xi}] f
+ {\cal N}[f].
\label{pdeqn7}
\end{eqnarray}
It should be mentioned here that this analysis already provides
an explicit estimate for the growth/shrinkage of amplitude and width
over time, given that we assume that $A_{\tau}/A=G=$const.
In particular, $A_t=A_{\tau} \tau_t= A_{\tau} A^{s-1}=G A^s$, which in
accordance to the considerations of the previous section leads
to the evolution of $A \sim (t^{\star}-t)^{1/(-s+1)}$. A similar
analysis can be performed for $L$ such that
$L_t=L_{\tau} \tau_t=L_{\tau} L^{-a}= -G \frac{s-1}{a} L^{1-a}$,
leading to $L \sim (t^{\star}-t)^\frac{1}{a}$.
As a result of this analysis, our pulse-like entity
touching (and potentially) crossing infinity will do
so in a self-similar manner. 
%The scaling evolution
%of the ``tip of the parabola'' leading to the crossing
%is briefly analyzed in a separate appendix in some
%generality.
%Before we embark into the detailed consideration of the computational
%framework, we explore the self-similar (MN) formulation of the relevant
%dynamics. Using 

For the specific example of Eq.~(\ref{pdeqn10}), 
using $w=A f(\xi,\tau)$,
we obtain that
\begin{eqnarray}
{\cal L}[\partial_{\xi}] w = \frac{A}{L^2} (w_{\xi \xi}-\frac{2}{w} w_{\xi}^2);
\quad {\cal N}[w]= A w + A^2 w^2
\label{pdeqn11}
\end{eqnarray}
It is then evident that the dynamics is not directly self-similar
(due to the different scaling of the two terms within ${\cal N}$),
but
%, similarly to the case explored in the Appendix II, 
only
{\em asymptotically self-similar}. When $w$ (and $A$) is small,
the exponential growth associated with the linear term is dominant.
However, as the amplitude increases, eventually the quadratic term
takes over, leaving the linear term as one of ever-decreasing-significance
``offending" to the exact self-similar evolution. When the linear
term becomes negligible, the self-similar evolution requires
that $A/L^2 = A^2$, providing the scaling of $A \sim 1/L^2$, i.e.,
in this case $s=2$ and $a=2$ for the general formulation above.
From there, all the scalings associated with self-similarity can
be directly deduced as explained previously.

%\clearpage

\section*{Supporting Information}

\subsection*{The linear ODE case}

%The main text of the paper must stand on its own without the SI. Refer to SI in the manuscript at an appropriate point in the text. Number supporting figures and tables starting with S1, S2, etc. Authors are limited to no more than 10 SI files, not including movie files. Authors who place detailed materials and methods in SI must provide sufficient detail in the main text methods to enable a reader to follow the logic of the procedures and results and also must reference the online methods. If a paper is fundamentally a study of a new method or technique, then the methods must be described completely in the main text. Because PNAS edits SI and composes it into a single PDF, authors must provide the following file formats only.

The mapping of the dynamics onto a circle can also be performed
for the case of the simple exponential (rather than the power law self-similar,
finite-time collapse) arising from the simple linear ODE of the form:
\begin{eqnarray}
\dot{x}= \pm x.
\label{eqnap1}
\end{eqnarray}
with the standard solution
\begin{eqnarray}
x(t)=e^{\pm (t-t^{\star})}.
\label{eqnap2}
\end{eqnarray}
Here the dynamics can be written in hyperbolic form as
\begin{eqnarray}
\left(\frac{e^{\mp (t-t^{\star}) } + x}{2}\right)^2
- \left(\frac{e^{\mp (t-t^{\star}) } - x}{2}\right)^2 =1,
\label{eqnap3}
\end{eqnarray}
and the variables
\begin{eqnarray}
X &=& \cos(\theta)=\frac{e^{\pm (t^{\star}-t)}-x}{e^{\pm (t^{\star}-t)}+x}
\label{eqnap4}
\\
Y &=& \sin(\theta)=\frac{2}{e^{\pm (t^{\star}-t)}+x},
\label{eqnap5}
\end{eqnarray}
can be defined so that $X^2+Y^2=1$. In fact, substituting the
exact solution of Eq.~(\ref{eqnap2}), it is straightforward
to realize that $X=\tanh(\mp (t-t^{\star}))$ and
$Y={\rm sech}(t-t^{\star})$, resulting in the circular dynamics
being a realization of the simple identity
$\tanh^2+{\rm sech}^2=1$.
%It might be interesting to take the view that such transformations 
%(including the case below), that map scalar ODEs to constant angular velocity dynamics, can %be thought of as transformations
%to curvature driven flow, here on a constant curvature manifold (the circle). 

%%\subsubsection*{SI Text}

%%Supply Word, RTF, or LaTeX files (LaTeX files must be accompanied by a PDF with the same file name for visual reference).

\subsection*{An asymptotically self-similar ODE case}
%{\bf An asymptotically self-similar ODE case}

We so far focused on genuinely self-similar examples; the
corresponding ideas can also be extended to {\em asymptotically} self-similar cases [1] %~\cite{barenb}
that are not genuinely self-similar
in that they possess ``offending'' terms, yet upon approaching
the singularity  the self-similar terms dominate,
with the offending ones playing a progressively less important role.
Our approach can easily be adapted to this case. \\

Our simple example variant here will
be of the form:
\begin{eqnarray}
\dot{x}=2 x + x^2.
\label{eqnap6}
% 6-->21
\end{eqnarray}
Direct integration again can yield the exact solution in the form:
\begin{eqnarray}
x(t)=\frac{2 e^{2 (t-t^{\star})}}{1-e^{2 (t-t^{\star})}}.
\label{eqnap7}
\end{eqnarray}
It can be seen (when integrating Eq.~(\ref{eqnap6}))
that in this case the observable $\log(x/(x+2)$
is the one that linearly crosses through $0$ (as $2 (t-t^{\star})$).
%
%Assuming
%that $x$ is 
For $x$ large, this quantity becomes
\begin{eqnarray}
\log(\frac{1}{1+\frac{2}{x}}) \approx -2 \frac{1}{x} + 2 \frac{1}{x^2}
-\frac{8}{3} \frac{1}{x^3}
\label{eqnap8}
\end{eqnarray}
Hence, indeed at large times, it is the quadratic term that takes
over since the dominant behavior of $x(t)$ is like $1/(t^{\star}-t)$.
However, as $t \rightarrow t^{\star}$, the relevant asymptotics reads:
\begin{eqnarray}
x(t) = \frac{1}{t^{\star}-t} -1 + \frac{t^{\star}-t}{3} - \frac{(t^{\star}-t)^3}{45} + \dots
\label{eqnap9}
\end{eqnarray}
enabling one to observe the explicit (lower order) contribution of
the terms offending to the self-similarity. 
the collapse time, denoted by $t^{\star}$, is still
determined by the initial data as $t^{\star}=-(1/2) \log(x(0)/(x(0)+2))$.

Nevertheless, in this case as well, our computational prescription
can be carried out. Eq.~(\ref{eqnap6}) can be integrated until
$x$ becomes large. we then revert to $y=1/x$ which has the straightforward
ODE dynamics:
\begin{eqnarray}
\frac{dy}{dt}=-2 y -1
\label{eqnap10}
\end{eqnarray}
(using the transformation to obtain the initial condition $y(0)$)
and the equally simple solution $y(t)=-1/2 + (y(0)+1/2) e^{-2 t}$.
The solution of the latter problem of Eq~(\ref{eqnap10})
crosses $0$ en route to its approach of the asymptotic value of
$-1/2$. Finally, once the infinity has been bypassed, we return to the
simulation of Eq.~(\ref{eqnap6}), as before.

{\it Mapping the dynamics to a circle}. The solution of Eq.~(\ref{eqnap6}) 
can be rewritten as:
\begin{eqnarray}
x(t)=\frac{2}{e^{2 (t^*-t)}-1} \Rightarrow x \frac{(e^{2 (t^*-t)}-1)}{2} =1.
\label{si1}
\end{eqnarray}
Using the compactification the exact same way as Eqs. [3] and [4] of
the main text 
and only replacing $t^*-t$ with:  $(e^{2 (t^*-t)}-1)/2$, the compactification scheme carries through.

%%%YGK
%%%  Costa, the Eqs (3) and (4) above are the (3) and (4) of the main text -
%%%  so you have to fix them to give the right numbers...
%%%   also make sure what the ???-??? equations are below

a) In this case, if $t \rightarrow t^*$, we Taylor expand and retrieve (from the first term) the limit of exactly Eqs. [3]--[4].
This is the contribution that stems from the $x^2$ term in the ODE.

b) In the case of $t \rightarrow 0$ (or anyway far from $t^*$) the exponential dominates and the (-1) coming
from the $x^2$ term is irrelevant. This is the contribution that stems from 
the $2x$ term in the ODE.

\subsection*{Time to transition}

Following numerous works
including [26, 27],%~\cite{sntime,bender}, 
we consider a radial contour along
the complex plane i.e., the arc of a circle from the real to the positive
imaginary axis. Then, along this arc (denoted by C),we have for  $T$, the elapsed time:
\begin{eqnarray}
T=\int_C dt = \int_C \frac{dz}{z^3} = \int_0^{\pi/2} \frac{R i e^{i \phi}}{R^3
e^{3 i \phi}} d \phi.
\label{eqn11}
% 11-->31
\end{eqnarray}
Bearing in mind the radial nature of the contour (which renders
$R$ constant), factoring out
$1/R^2$ and taking the limit as $R \rightarrow \infty$, we obtain
a vanishing result, even though the angular integral amounts
to unity. I.e., interestingly, it takes a finite time to reach
from everywhere along the real axis an equidistant point along the
imaginary axis, yet this time vanishes as we approach infinity,
in line with the analytical result for $\dot{x}=x^3$. 
In the case of $\dot{z}=z^2$, there is a similar
result justifying the infinitesimal time of return there
from the positive to the negative real axis.

\section*{Acknowledgments}
The authors gratefully acknowledge support from the US NSF
and the US Air Force Office of Scientific Research (Dr. F. Darema), 
as well as stimulating discussions with Profs. D. Aronson, C. Bender, E. Bollt, P. Constantin, M. Dafermos, C. Tully and Z. Musslimani.

\nolinenumbers

%This is where your bibliography is generated. Make sure that your .bib file is actually called library.bib
%\bibliography{library}

\end{document}